\documentclass{gtart}
\usepackage{amssymb,rlepsf}
\let\relabela\adjustrelabel

\input gtoutput
\volumenumber{2}\papernumber{9}\volumeyear{1998}
\pagenumbers{175}{220}\published{4 January 1999}
\proposed{David Gabai}
\seconded{Wolfgang Metzler, Cameron Gordon}
\received{3 July 1997}\revised{9 January 1998}
\accepted{4 January 1999}

\def\reals{\mathbb{R}}

\def\ints{\mathbb{Z}}

\def\eps{\epsilon}

\newtheorem{theorem}{Theorem}[section]
\newtheorem{hthm}[theorem]{Halting Theorem}
\newtheorem{prop}[theorem]{Proposition}
\newtheorem{lemma}[theorem]{Lemma}
\newtheorem{cor}[theorem]{Corollary}
\newtheorem{remark}[theorem]{Remark}
\let\pf\proof
\let\endpf\endproof
\newcommand{\med}{\medskip}

\def\makemath#1{\ifmmode{#1}\else${#1}$\fi}

\def\cD{\makemath{\cal D}}

\def\cF{\makemath{\cal F}}
\def\cH{\makemath{\cal H}}

\def\cW{\makemath{\cal W}}

\begin{document}

\title{A new algorithm for recognizing the unknot}%

\author{Joan S Birman\\Michael D Hirsch}

\address{Math Dept, Columbia University, NY, NY 10027, USA\\
Math and CS, Emory University, Atlanta, GA 30322, USA}

\email{jb@math.columbia.edu\\hirsch@mathcs.emory.edu} 

\begin{abstract}

  The topological underpinnings are presented for a new algorithm
  which answers the question: ``Is a given knot the unknot?''  The
  algorithm uses the braid foliation technology of Bennequin and of
  Birman and Menasco.  The approach is to consider the knot as a
  closed braid, and to use the fact that a knot is unknotted if and
  only if it is the boundary of a disc with a combinatorial foliation.
  The main problems which are solved in this paper are: how to
  systematically enumerate combinatorial braid foliations of a disc;
  how to verify whether a combinatorial foliation can be realized by
  an embedded disc; how to find a word in the the braid group whose
  conjugacy class represents the boundary of the embedded disc; how to
  check whether the given knot is isotopic to one of the enumerated
  examples; and finally, how to know when we can stop checking and be
  sure that our example is {\em not} the unknot.

\end{abstract}

\asciiabstract{%
The topological underpinnings are presented for a new algorithm
which answers the question: `Is a given knot the unknot?'  The
algorithm uses the braid foliation technology of Bennequin and of
Birman and Menasco.  The approach is to consider the knot as a
closed braid, and to use the fact that a knot is unknotted if and
only if it is the boundary of a disc with a combinatorial foliation.
The main problems which are solved in this paper are: how to
systematically enumerate combinatorial braid foliations of a disc;
how to verify whether a combinatorial foliation can be realized by
an embedded disc; how to find a word in the the braid group whose
conjugacy class represents the boundary of the embedded disc; how to
check whether the given knot is isotopic to one of the enumerated
examples; and finally, how to know when we can stop checking and be
sure that our example is not the unknot.}




\maketitlepage

\section{Introduction}
\label{section:Introduction}
The goal of this manuscript is the development of a new algorithm to
answer the question: Given a knot $K$ which is defined by a diagram,
does $K$ represent the unknot? Our algorithm is suitable for computer
enumeration.  Its approach is straightforward:
\begin{enumerate} 
\item We show how to construct a sufficiently large set of diagrams
  which represent the unknot;
\item We introduce a complexity function which allows us to order
  these diagrams, as we construct them, in order of complexity;
\item We learn how to test in a systematic way whether an arbitrary
  diagram for a knot $K$ is equivalent to one of the diagrams on the
  list;
\item We arrange that the checking process stop in a finite time. 
\end{enumerate}
The focus of this paper will be on the topological underpinnings of
the algorithm. We are in the process of implementing the algorithm,
and of assembling computer-generated data for the ordered list which
we produce, and the data should be interesting. We plan to write a
second paper on that work, when it is complete and done efficiently
enough to give us the data we would like to see.

Before we describe our approach, we give a brief review of related
work on the problem.

\begin{itemize}
\item In the 1960's W Haken \cite{Haken} used the concept of a normal
  surface to show that the homeomorphism problem is solvable for
  triangulated 3--manifolds which contain an incompressible surface.
  The unknot $K$ bounds a disc, and the disc is an incompressible
  surface in the 3--manifold which is obtained by
  removing an open tubular neighborhood of $K$ from $S^3$, so Haken's
  work proves the existence of an algorithm for recognizing the
  unknot. See \cite{Hass} for a review of Haken's contribution to the problem.
\item Several months after an earlier draft of this paper was
  submitted for publication Hass, Lagarias and Pippenger announced new
  results on the problem of recognizing the unknot in \cite{HLP}.
  Their work investigates Haken's algorithm in detail, using more
  recent contributions of Jaco and Oertel and new techniques from
  computer science to sharpen Haken's algorithm and place explicit
  bounds on its running time. They then go on to prove that the
  problem of recognizing the unknot is in class NP.
That issue is somewhat tangential to the subject matter of this paper.
  The emphasis of this paper is on the existence of a new algorithm
  for unknottedness, rather than the complexity.  
  Our proof that we can stop checking after an explicit finite time
  will use similar methods to \cite{HLP}, but the algorithm itself is
  totally different.  
It is not clear at this writing whether our algorithm
demonstrates that the problem of detecting unknots is in NP.
  
\item In a different direction, an open conjecture is whether there
exists a non-trivial knot whose Jones polynomial is 1. If none
exists then the Jones polynomial detects the unknot. 

\end{itemize} 
Our work is in the setting of closed braids.  When we say that we list
`a sufficiently large set of diagrams which represent the unknot' we
mean that we enumerate the conjugacy classes of the unknot in the
braid group $B_n$, where $n$ is the number of Seifert circles in the
given diagram of $K$, in an appropriate order.  Since there is a very simple
algorithm \cite{Vogel} to change every knot diagram with
$n$ Seifert circles to a closed $n$--braid diagram and to read off a
representing open braid, and since the conjugacy problem in the braid
group is a solved problem (\cite{Garside},\cite{ElM},\cite{BKL}) we
will then have the tool we need to solve problem 3, ie to test (one
conjugacy class at a time) whether the closed $n$--braid which we
constructed from the given diagram of $K$ is conjugate to one of the members of
our list. The list is infinite for $n\geq 4$, and our main problems are to
construct the list, to order it in an appropriate way, and to learn
when we can stop testing. Those are non-trivial problems, and we will
bring much structure to bear on them, in addition to using the known
solution to the conjugacy problem.

The unknot is the unique knot which bounds a disc, and our tool for
enumerating its closed $n$--braid representatives is based on the
combinatorics of certain {\em braid foliations} of the disc.  These
foliations were introduced by D Bennequin in \cite{Bennequin}. They
were studied systematically as a tool in knot theory by the first
author and W Menasco in a series of papers with the common title
`Studying links via closed braids', for example see
\cite{BMV}. Our detailed work involves many ideas from those papers, but for
convenience our references will be mainly to the review article \cite{BF}, which
gathers together in one place the machinery developed in those papers. Our
technique for enumerating all closed braid representatives of the unknot is
in fact implicit in D. Bennequin's work
\cite{Bennequin}.  It is a method of ``stabilizing'' a complicated 
embedded disc
to obtain a simpler one whose boundary has much higher braid index.  
We use the
reverse of this procedure to generate complicated discs whose boundaries 
have low
braid index. Some of these are not embeddable, so we develop a new
(and surprisingly simple) test for embeddability, allowing us to
eliminate any non-embeddable ones that may have arise in the course of
the enumeration.

Having in hand an embeddable foliated disc with associated
combinatorial data, we know (from a theorem proved in
\cite{BF}) that its embedding in 3--space relative to cylindrical
coordinates is determined up to foliation-preserving isotopy.
However, we still have to solve the problem of determining a word in
the generators of the braid group which represents the boundary of the
given disc.  The solution to that problem is new to this paper and
turns out to be quite elegant.  All of this, combined with the
solution to the conjugacy problem in $B_n$ in \cite{BKL}, solves
problems 1, 2 and 3 above.
  
%%% Replaces page 4, line -5 to page 5, line 8 
The complexity measure which we assign to our foliated disc is a pair
of integers $(n,v)$, where $n$ is the braid index of the boundary and
$v$ is the number of times the braid axis intersects the disc.  The
given knot $K$ is defined by a diagram, and as noted earlier $n$ is also
the number of Seifert circles in the diagram, which is thus fixed for
each example.  To find a bound on $v$ we must relate $v$ to the crossing number
$k$ of $K$. For this part of the argument we construct a triangulation of the
complement  of a tubular neighborhood of $K$, doing it so that the braid axis
meets the interiors of exactly 4 tetrahedra, in a controlled way. We then
find an upper bound on the number $t$ of tetrahedra, as a function of $k$ and
$n$.  As is well-known, the work of Kneser and Haken
implies the existence of an upper bound on the number of times the disc we are
seeking can intersect a single tetrahedron. Fortunately we do not need to
compute that bound because  Lemma 6.1 of \cite{HLP} does the job for us.
Multiplying by 4 we obtain an upper bound for $v$, which then tells us when we
can stop testing.
%%% end of modification.  Old version is below.

% The complexity measure which we assign to our foliated disc is a pair
% of integers $(n,v)$, where $n$ is the braid index of the boundary and
% $v$ is the number of times the braid axis intersects the disc.  The
% given knot $K$ is defined by a diagram, and as noted earlier $n$ is
% the number of Seifert circles in the diagram, which is thus fixed for
% each example.  To place a bound on $v$ we must relate $v$ to
% the crossing number of $K$. For this part of the argument we are
% indebted to ideas of Hass, Lagarias and Pippenger \cite{HLP}. As in
% \cite{HLP}, we construct a triangulation of the complement of a tubular
% neighborhood of $K$, but we do it so that the braid axis meets the
% interiors of only 4 tetrahedra, and does so in a controlled way. The
% maximal number of tetrahedra in the triangulation is bounded above by
% a function which depends upon the crossing number of the given diagram
% for $K$ and $n$. Knowing this bound, a theorem in \cite{HLP} gives us
% a bound for the number of times the disc which we are seeking can
% intersect the 4 distinguished tetrahedra. From this we obtain an upper
% bound for $v$, which tells us when we can stop testing.

Here is an outline of the paper. In Section \ref{section:braid foliations
  of spanning surfaces for knots} we review the prerequisite material
about braid foliations. 
This section is without proofs, but the
material is fairly understandable and believable.  A convenient
reference, complete with proofs, is available \cite{BF}.  In
Section \ref{section:testing for embeddability} we show how to test
whether a given combinatorially foliated surface actually corresponds
to an embedded foliated surface, and if so how to find a braid word
which represents the boundary. In Section \ref{section:enumerating} we
show how to enumerate the ordered list of closed $n$--braid
representatives of the unknot which is the basis for our algorithm.
In Section \ref{section:halting} we prove our `halting theorem'.  In
Section \ref{section:the algorithm} we present the algorithm.

For completeness, we show in Appendix A how to rapidly modify an
arbitrary knot diagram to a closed braid diagram, with control over
the extra crossings which are added. This part of the algorithm is
based upon the work of Vogel, reported on in \cite{Vogel}.  In
Appendix B we review the solution to the conjugacy problem in $B_n$
which we are using in our algorithm. The theoretical basis for that
algorithm is established in \cite{BKL}.

\med
{\bf Acknowledgements}\qua We thank Elizabeth Finkelstein for her many
contributions to this paper, both through her work in \cite{BF} and
through our discussions with her at an early stage in this work. We
thank Joel Hass for helpful conversations. It was only after we read
the manuscript \cite{HLP} and talked to him that we understood how to
give the proof of the Halting Theorem which is
presented here, replacing a much more awkward solution to that problem
in the earlier draft of this paper. We also thank William Menasco and
Brian Mangum for helpful conversations.

The first author acknowledges partial support from the following
sources: the US National Science Foundation, grants DMS 94-02988 and DMS
97-05019; Barnard College, for salary support during a Senior Faculty Research
Leave; the Mathematical Sciences Research Institute, where she was a
Visiting Member when part of this work was done; and the US Israel
Binational Science Foundation. The second author would like to thank
the US National Science Foundations for partial support under grants
ASC-9527186 and DMS-9404261.

\section{Braid foliations of spanning surfaces for knots} 
\label{section:braid foliations of spanning surfaces for knots}

In this section we will review the basic theorems about the braid
foliations associated to an  incompressible orientable surface
spanning a knot or link.  The main theorems in this field are collected in
the survey paper \cite{BF}, which we will use as our preferred reference for
the material of this section.

While we are interested mainly in closed braid representatives of the
unknot, we state everything in terms of links, and the (possibly
disconnected) incompressible spanning surfaces of genus $g\geq 0$
which they bound.  Admitting these complications requires almost no
extra work and yields much simplified proofs in
Section~\ref{section:testing for embeddability}.  It is also necessary
for any future extensions of these ideas to true (knotted) knots.

Choose cylindrical coordinates $(r,\theta,z)$ in 3--space ($S^3$
thought of as $\reals^3$ union a point at $\infty$). A link $K$ is a
{\em closed braid} with the $z$--axis as its {\em braid axis} if each
component of $K$
has a parameterization $\{(r(t),\theta (t), z(t)) : t\in [0,1]\}$ such
that $r(t)\neq 0$ and $d\theta /dt >0$ for all $t\in [0,1] $. This
means that $K$ intersects each half-plane $H_\theta$ through the axis
transversally, as in Figure \ref{figure:closed
braid}.  It follows that the number of points in $K \cap H_\theta$ is
independent of $\theta$. This number is the {\em braid index}, $n =
n(K)$, of $K$. Notice that the closed 4--braid diagram in Figure 1 contains 4 Seifert
circles. See the Appendix for a proof that any diagram with $n$ Seifert circles can be
modified to a closed $n$--braid, adding a controlled number of crossings. 

\begin{figure}[htb!]
\centerline{\relabelbox\small
\epsfxsize 3 in\epsfbox{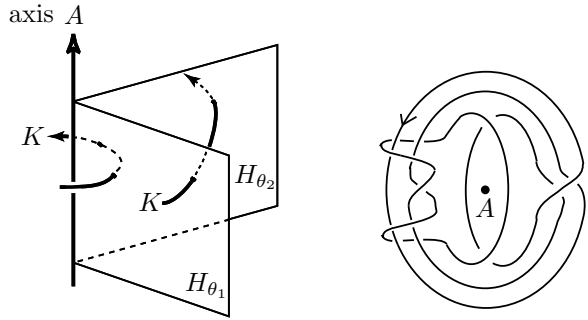}
\relabel {A}{$A$}
\relabela <-2pt, -2pt> {A1}{$A$}
\relabela <-4pt, 0pt> {K}{$K$}
\relabela <-4pt, 1pt> {K1}{$K$}
\relabela <-2pt, 1pt> {axis}{axis}
\relabela <-4pt, 0pt> {H}{$H_{\theta_1}$}
\relabela <-2pt, 0pt> {H2}{$H_{\theta_2}$}
\endrelabelbox}
\caption{Braids and closed braids
\label{figure:closed braid}}
\end{figure}

Let $K$ be a closed braid and let $F$ be an embedded orientable
surface of minimal genus spanned by $K$.  Note that both $F$ and
$\partial F$ may be disconnected.  Let $\cH = \{H_\theta \ : \ 
\theta\in [0,2\pi ]\}$ be the open book decomposition of $\reals^3$ by
half-planes with boundary on the $z$--axis. Assume that $F$ is in
general position with respect to $\cH$ and consider the induced
foliation on $F$. This foliation is the {\em braid
  foliation} of $F$.

The braid foliation is singular both at points where
$F$ meets the braid axis $A$, and at points where $F$ is tangent to leaves of
$\cH$. By general position, these latter type of singularities can be assumed, 
{\it a priori}, to be of saddle
type, or center type with neighborhoods foliated by circles. The following
is a restatement of results of Bennequin \cite{Bennequin}. It classifies the
leaves and singularities of  the braid foliation after an isotopy.  

\begin{theorem} \label{theorem:classification}
{\rm (\cite{BF}, Theorems 1.1,1.2) } \ After a modification of $F$ rel
$\partial F$ (by isotopy when $\partial F$ is non-split) the following hold:
\begin{enumerate}
\item The braid foliation  near ${\partial F}$ (see Figure
\ref{figure:transverse near the boundary}) is transverse to
$\partial F$. 
It is radial in a neighborhood of each point of $A\cap F$ (see 
Figure~\ref{figure:radial near the tile vertices}).

\begin{figure}[htb!]
\begin{minipage}{0.02\textwidth}
\hfill
\end{minipage}
\begin{minipage}{0.45\textwidth}
\centerline{\relabelbox\small
\epsfxsize 2 in\epsfbox{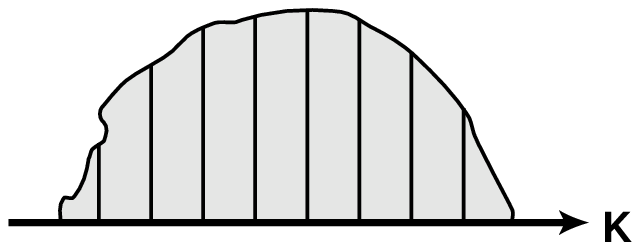}
\relabela <0pt, 1pt> {K}{$K$}
\endrelabelbox}\vglue -3mm
\caption{The braid foliation in a neighborhood of $\partial F$.
\label{figure:transverse near the boundary}}
\end{minipage}
\begin{minipage}{0.03\textwidth}
\hfill
\end{minipage}
\begin{minipage}{0.45\textwidth}\vglue 5mm
\centerline{\relabelbox\small
\epsfxsize 2 in\epsfbox{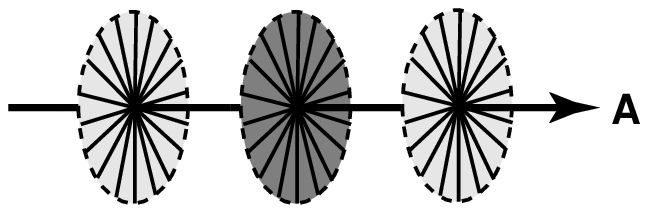}
\relabela <0pt, 0.5pt> {A}{$A$}
\endrelabelbox}
\caption{The braid foliation in a neighborhood on $F$ of each point of $A\cap F$.
\label{figure:radial near the tile vertices}}
\end{minipage}
\end{figure}

\item  \label{item:nonsingular}
The non-singular leaves of the braid foliation fall into three types: $a$,
$b$ and $c$.  (See
the left sketch in Figure~\ref{figure:arctypes}.)
\begin{itemize}
\item [$a$:] arcs with one endpoint on $K$ and one on $A$, and
\item [$b$:] arcs with both endpoints on $A$.
\item [$c$:] arcs with both endpoints on $K$.
\end{itemize}
However arcs of type $c$ are necessarily singular because
intersections of $K$ with every fiber of $\cH$ are coherently oriented
(see the right sketch in Figure~\ref{figure:arctypes}), so they do not
occur as non-singular leaves.

\begin{figure}[htb!]
\centerline{\relabelbox\small
\epsfxsize 4 in\epsfbox{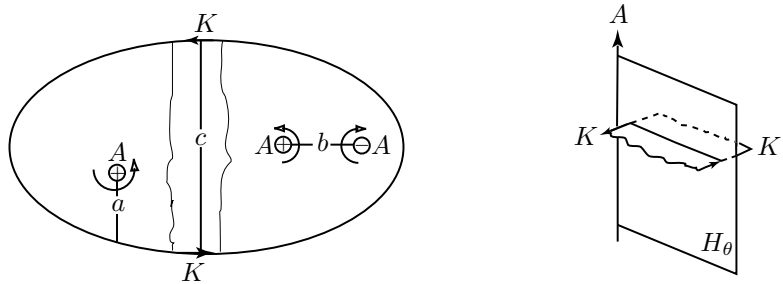}
\relabela <2pt, -2pt> {K1}{$K$}
\relabela <-6pt, 0pt> {K2}{$K$}
\relabela <0pt, 0pt> {K3}{$K$}
\relabela <-3pt, 0pt> {K4}{$K$}
\relabela <-3.5pt, 0.5pt> {H}{$H_\theta$}
\relabela <-2pt, 0pt> {A1}{$A$}
\relabela <-1.5pt, 0pt> {A2}{$A$}
\relabela <-1pt, 0pt> {A3}{$A$}
\relabela <0pt, 1pt> {A4}{$A$}
\relabela <-1.8pt, 1.6pt> {m}{$\scriptstyle -$}
\relabela <-1.85pt, 0pt> {p}{$\scriptstyle +$}
\relabela <-1.8pt, 0.5pt> {p1}{$\scriptstyle +$}
\relabela <-0.5pt, 0.5pt> {a}{$a$}
\relabela <-0.3pt, 0pt> {b}{$b$}
\relabela <0pt, 0pt> {c}{$c$}
\endrelabelbox}\vglue -9mm
\caption{Leaves of type $a$,$b$ and $c$
\label{figure:arctypes}}
\end{figure}
   
\item The singular points of the braid foliation are of two types which we call
\underline{vertices} and \underline{singularities}: 
\begin{itemize}
\item [] 
\hspace{-.25 truein}
Vertices: points of intersection between $F$ and $A$.  The
foliation is radial  near the vertices 
(see Figure~\ref{figure:radial near the tile vertices}). 

\item [] \label{item:singular}
\hspace{-.25 truein}
Singularities: A point of tangency between $F $ and some
$H_\theta $.  These tangencies are simple saddles
(see Figure~\ref{figure:3 tile types}).
Such $H_\theta $ are called  singular.  The point of tangency is
called a singularity.  The singularity, together
with its four leaves, is called a \underline{singular leaf}.  
The four leaves of a
singular leaf are called \underline{branches} of the singular leaf.
\end{itemize}

%%% mdh moved figure so we don't have only two lines of text on page

\item Singularities fall into three types: $aa$, $ab$ and $bb$ (see
Figure~\ref{figure:3 tile types}).
\begin{enumerate}
\item [$aa$:] those singularities between  two $a$--arcs, 
\item [$ab$:] those between an $a$--arc and a $b$--arc, and 
\item [$bb$:] those between two $b$--arcs. 
\end{enumerate}

%%%mdh figure moved to here.

\begin{figure}[htb!]
\centerline{\relabelbox\small
\epsfxsize 4 in\epsfbox{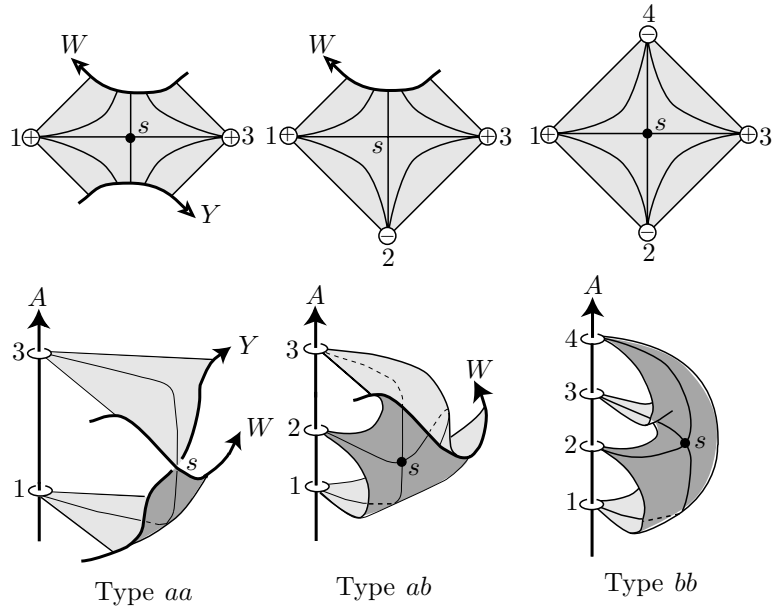}
\relabela <-2pt, 0pt> {A1}{$A$}
\relabela <-2pt, 0pt> {A2}{$A$}
\relabela <-2pt, 0pt> {A3}{$A$}
\relabela <0pt, 0pt> {Type aa}{Type $aa$}
\relabela <0pt, 0pt> {Type ab}{Type $ab$}
\relabela <0pt, 0pt> {Type bb}{Type $bb$}
\relabela <0pt, 0pt> {s1}{$s$}
\relabela <0pt, 0pt> {s2}{$s$}
\relabela <0pt, 0pt> {s3}{$s$}
\relabela <-5pt, -4pt> {s4}{$s$}
\relabela <0pt, 0pt> {s6}{$s$}
\relabela <0pt, 0pt> {s7}{$s$}
\relabela <0pt, 0pt> {W1}{$W$}
\relabela <0pt, 2pt> {W2}{$W$}
\relabela <0pt, 2pt> {W3}{$W$}
\relabela <0pt, 0pt> {W4}{$W$}
\relabela <0pt, 0pt> {Y1}{$Y$}
\relabela <0pt, 0pt> {Y2}{$Y$}
\relabela <0pt, 0pt>  {11}{1}
\relabela <0pt, -1pt>  {12}{1}
\relabela <0pt, -1pt>  {13}{1}
\relabela <-1pt, 0pt>  {14}{1}
\relabela <0pt, -1pt>  {15}{1}
\relabela <-1pt, 0pt>  {16}{1}
\relabela <0pt, 0pt>  {21}{2}
\relabela <0pt, -2pt>  {22}{2}
\relabela <0pt, -2pt>  {23}{2}
\relabela <-2pt, 0pt>  {24}{2}
\relabela <0pt, 0pt>  {31}{3}
\relabela <-1pt, -1pt>  {32}{3}
\relabela <-1pt, 0pt>  {33}{3}
\relabela <0pt, -1pt>  {34}{3}
\relabela <0pt, -1pt>  {35}{3}
\relabela <-1pt, 0pt>  {36}{3}
\relabela <0pt, 0pt>  {41}{4}
\relabela <0pt, 0pt>  {42}{4}
\relabela <-2pt, 2pt> {m3}{$\scriptstyle-$}
\relabela <-2pt, 2pt> {m4}{$\scriptstyle-$}
\relabela <-2pt, 2pt> {m5}{$\scriptstyle-$}
\relabela <-1pt, 1pt> {p2}{$\scriptstyle+$}
\relabela <-1pt, 1pt> {p3}{$\scriptstyle+$}
\relabela <-1pt, 1pt> {p4}{$\scriptstyle+$}
\relabela <-1pt, 1pt> {p5}{$\scriptstyle+$}
\relabela <-1pt, 1pt> {p6}{$\scriptstyle+$}
\relabela <-1pt, 1pt> {p7}{$\scriptstyle+$}
\endrelabelbox}
\caption{The three singularity types
\label{figure:3 tile types}}
\end{figure}

\item \label{item:orderings}
The vertices are (circularly) ordered by their order on the
braid axis $A$.  After isotopy
distinct singular leave are on distinct $H_\theta $ and are thus
circularly ordered, as well.
\end{enumerate}
\end{theorem}

 In part \ref{item:nonsingular} of the theorem, a salient point is
that non-singular circular leaves can occur at local minima and
maxima, and the theorem says that (subject to the assumption
that $F$ has minimum genus among all orientable surfaces bounded by
$K$) we can cut off any maxima and fill
up the minima without leaving the class of surfaces which are of interest.  If
there no local minima or maxima, only the vertex and saddle type singularities
of part \ref{item:singular} can occur.  This arrangement, together with the fact
that a leaf with both of its endpoints on $K$ is necessarily singular, 
%%% mdh  changed "are" to "is"
is
responsible for the rich combinatorics of braid foliations.    

There are additional combinatorial data.  Since the
original braid $K$ was oriented, $F$ has an orientation.  If that
orientation agrees with the orientation of $A$ at a vertex $v$,
(ie the normal vector has positive inner product with the oriented
tangent vector to $A$) then we say the vertex is {\em positive}, otherwise
it is {\em negative}.  Similarly, at a singularity, if the normal vector is
positive with respect to $d\theta $ we say the singularity is
{\em positive}, otherwise it is {\em negative}.

There is dual viewpoint with regard to these foliations.  Instead of focusing
on singular leaves, one can break the surface up
into {\em tiles}, one for each singularity, by cutting along appropriate
$a$--arcs and
$b$--arcs. Figure~\ref{figure:3 tile types} shows the tile types, for each
singularity.  This is the point of view
in \cite{BF}, and so we use the terminology ``tiled surface''.

We will need the following facts about singularities. They are self
evident from Figure~\ref{figure:3 tile types} and follow easily from
orientation considerations and the previous theorem.  Each singularity
is connected to exactly two positive vertices along non-adjacent
branches of the singular leaves.  Each of the other two leaves can go
to either the boundary of the orientable surface, or to a negative
vertex.  
%%% mdh We don't need this anymore.
%If we
%consider a point of intersection of a branch $b$ of a singular leaf
%and the braid, then there is a well-defined branch which is in the
%clockwise direction from $b$ when viewed from the positive side of
%$F$, and another which is counterclockwise from $b$.  For instance,
%in Figure~\ref{figure:3 tile types}, consider the point where the
%singular leaf through the point  $s$ meets $W\subset K$ in the pictures of type
%$aa$ and type
%$ab$ singularities.  Call the branch of intersection $b$.  Vertex 1 is
%on the branch to the counterclockwise of $b$ and Vertex 3 is on the
%branch to the clockwise of $b$.

\med{\bf Definition}\qua A {\em tiled surface} $\cF$ is a 3--tuple $(F,G,C)$
where $F$ is an oriented surface, $G$ is a graph which is embedded in
$F$ (so that it has a well-defined neighborhood in F), with some
additional combinatorial data $C$ which we call {\em decorations}. The
graph should be of the type attainable as the graph of singular
leaves, ie, $G$ must be tripartite with each node either a {\em
  vertex}, a {\em singularity}, or a {\em boundary point}.  The
singularities are of index 4 and are adjacent to at most two boundary
points (on non-adjacent edges).  The boundary points are on the
boundary of $F$ and adjacent to a single singularity. Each component
$F' \subset F\setminus G$ is a disc, also $\partial F'$ contains
either exactly 4 edges of $G$ or exactly 3 edges, one of which meets
$F'$ on both sides.  Each vertex has a sign, as does each singularity.
The vertices are circularly ordered, as are the singularities.
Vertices of the same sign can be adjacent to the same singularity only
if they are on edges which are non-adjacent at the singularity.

We wish to emphasize the decorations of the tiled surface.  In the
tiled surface literature these decorations have been largely ignored,
or, at best, implicit.  Tiled surfaces were studied primarily in terms
of the graph with only secondary thought given to the decorations.

If \cF is a tiled surface, then $F$ can be foliated (uniquely, up to
homeomorphism) by a singular foliation so that $G$ is the union of
singular leaves.  Thus every tiled surface as defined above is
implicitly foliated by $a$--arcs and $b$--arcs as specified in
Theorem~\ref{theorem:classification}.  Traditionally, one
uses the foliation instead of the graph.  Graphs are more natural to use
in an algorithmic context, and make the decorations easier to
understand, so we use them in our definition.

It is natural to abuse notation and think of \cF as a surface, rather
than as a tuple $(F,G,C)$; we will try not to do this.  We shall
consistently use the convention that the symbol for the surface will
be a roman capital and the tiled surface will be the same letter in a
calligraphic font.

\med{\bf Definition}\qua An {\em embeddable tiled surface} is a tiled
surface which is actually achieved as the graph of singular leaves of
some embedded orientable surface with closed braid boundary.  This
embedding is essentially unique.  See Theorem~\ref{theorem:unique
embedding}.

One more concept from the previously cited papers of Bennequin, Birman
and Menasco will be helpful, because it gives a very simple way to
exclude unwanted examples. As before, we refer the reader to \cite{BF}
for a detailed exposition.  Consult Figure \ref{figure:essential
b-arcs}.

\begin{figure}[htb!]
\centerline{\relabelbox\small
\epsfxsize 4 in\epsfbox{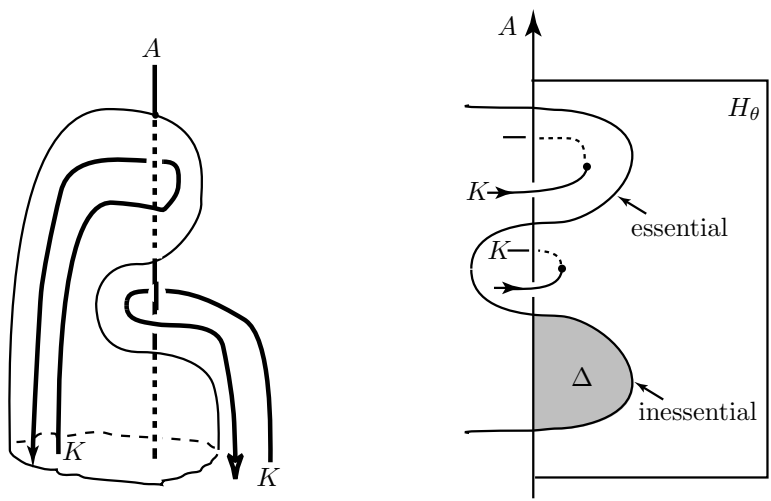}
\relabela <-2pt, 0pt> {A1}{$A$}
\relabela <-2pt, -9pt> {A2}{$A$}
\relabela <4pt, 0pt> {essential}{essential}
\relabela <2pt, 0pt> {inessential}{inessential}
\relabela <0pt, 0pt> {H}{$H_\theta$}
\relabela <0pt, 0pt> {D}{$\Delta$}
\relabela <-2pt, 1pt> {K1}{$K$}
\relabela <-4pt, 1pt> {K2}{$K$}
\relabela <-2pt, 2pt> {K3}{$K$}
\relabela <-2pt, 1pt> {K4}{$K$}
\endrelabelbox}
\caption{Essential and inessential $b$--arcs
\label{figure:essential b-arcs}}
\end{figure}

The sketch on the left illustrates a `pocket' in an embedded disc. It
cannot be removed because the knot is in the way. If the knot was not
an obstruction, we could eliminate the pocket (and remove two vertices
in the tiling) by an isotopy.  This leads us to a definition.

\med{\bf Definition}\qua A $b$--arc $\beta$ is said to be {\em essential} if
both sides of $H_\theta$ split along $\beta$ are pierced by $K$.  See
the right sketch in Figure \ref{figure:essential b-arcs}.  An embeddable
tiled surface is an {\em essential tiled surface} if all the $b$--arcs of the
braid foliation induced by the embedding of the tiled surface are
essential.  
An embeddable tiled surface which is not an essential tiled surface is
said to be 
an \emph{inessential tiled surface}.

\begin{figure}[htb!]
\centerline{\relabelbox\small
\epsfxsize 2.7 in\epsfbox{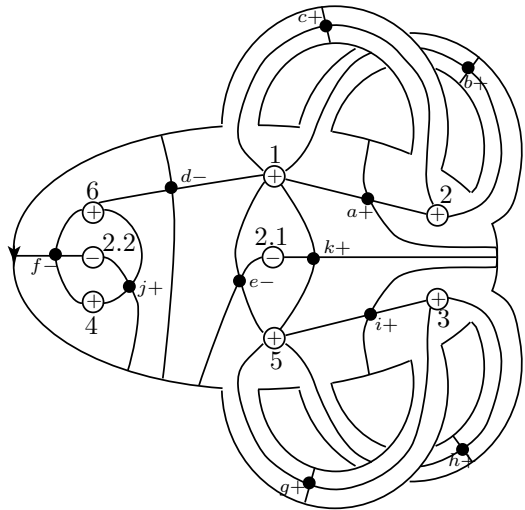}
\relabela <-1pt, 3.5pt> {m1}{$\scriptstyle-$}
\relabela <-0.6pt, 3.5pt> {m2}{$\scriptstyle-$}
\relabela <0pt, 1.7pt> {p1}{$\scriptstyle+$}
\relabela <0pt, 1.7pt> {p2}{$\scriptstyle+$}
\relabela <0pt, 1.7pt> {p3}{$\scriptstyle+$}
\relabela <-0.2pt, 1.5pt> {p4}{$\scriptstyle+$}
\relabela <0pt, 1.7pt> {p5}{$\scriptstyle+$}
\relabela <0pt, 1.7pt> {p6}{$\scriptstyle+$}
\relabela <0pt, 0pt>  {1}{1}
\relabela <0pt, 0pt>  {2}{2}
\relabela <0pt, 0pt>  {21}{2.1}
\relabela <0pt, 0pt>  {22}{2.2}
\relabela <0pt, 0pt>  {3}{3}
\relabela <0pt, 0pt>  {4}{4}
\relabela <0pt, 0pt>  {5}{5}
\relabela <0pt, 0pt>  {6}{6}
\relabela <-2pt, -1pt>  {a}{$\scriptstyle a+$}
\relabela <-0.5pt, 0.5pt>  {b}{$\scriptstyle b+$}
\relabela <-2pt, 1pt>  {c}{$\scriptstyle c+$}
\relabela <0pt, 0pt>  {d}{$\scriptstyle d-$}
\relabela <0pt, 0pt>  {e}{$\scriptstyle e-$}
\relabela <-3pt, 0pt>  {f}{$\scriptstyle f-$}
\relabela <-2.5pt, 0pt>  {g}{$\scriptstyle g+$}
\relabela <-0.5pt, 0pt>  {h}{$\scriptstyle h+$}
\relabela <0pt, 0pt>  {i}{$\scriptstyle i+$}
\relabela <0pt, 0pt>  {j}{$\scriptstyle j+$}
\relabela <0pt, 0pt>  {k}{$\scriptstyle k+$}
\endrelabelbox}
\caption {An example of an embeddable essential tiled surface
\label{figure:an example of an embeddable tiled surface}}
\end{figure}

\med{\bf Example}\qua Figure \ref{figure:an example of an embeddable tiled
  surface} gives an example of a tiled surface of genus 2. We will see
shortly how to test that it is embeddable and essential. The 11 singularities are
indicated by small black dots, signed and labelled
$a,b,c,d,e,f,g,h,i,j,k$ to correspond to their cyclic order in the
fibers around the axis.  The 8 white circles (6 positive and 2
negative) are the signed tile vertices. They are labelled
$1,2,2.1,2.2,3,4,5,6$ to them describe their order on the braid axis.
(It will become clear as we proceed why we choose non-integer labels
for the negative vertices).  

In the next section we will learn how to test whether a given example
is embeddable. We are aided in that project by the fact that when a
tiled surface is embeddable, then there is a unique
embedding, up to foliation preserving isotopy:

\begin{theorem} {\rm \cite[Theorem 4.1]{BF}}\qua 
\label{theorem:unique embedding}
The combinatorial data for a embeddable tiled surface $\cF$, ie, the
embedded graph $G$ and its embedding in $F$, the circular ordering
for its vertices, the circular ordering for its singularities, and the
signs of the vertices and singularities, determine the embedding in
$S^3$.  This embedding is unique up to foliation preserving isotopy.
The embedding of the boundary is determined by the same data,
restricted to singular leaves which meet the boundary and their
associated vertices.
\end{theorem}

\section{Testing for embeddability and finding the boundary word }
\label{section:testing for embeddability}
Given a knot or link, it is natural to ask what surfaces it spans.  In
this section we study a dual question: Given a tiled surface, is it
embeddable?  And if it is embeddable, what braid is represented by its
boundary?  Our main results on these matters are Theorems
\ref{theorem:boundary word} and \ref{theorem:embeddability test}. During most
of the section we will ignore the question of whether the surface is essential,
but at the end of the section  Proposition 
\ref{proposition:inessential b-arcs} will give a very simple test
which can be used to 
eliminate inessential tiled surfaces.

\subsection{A special case: positive tiled surfaces}
\label{subsection:the special case of positive tiled surfaces}
\begin{figure}[htb!]
\centerline{\relabelbox\small
\epsfxsize 4 in\epsfbox{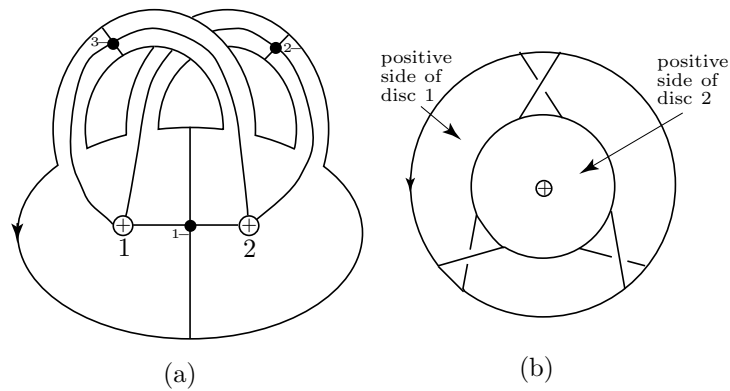}
\relabela <0pt, 0pt>  {1}{1}
\relabela <0pt, 0pt>  {2}{2}
\relabela <0pt, 0pt>  {a}{(a)}
\relabela <0pt, 4pt>  {b}{(b)}
\relabela <-0.5pt, 1.1pt> {p1}{$\scriptstyle+$}
\relabela <-0.5pt, 1.1pt> {p2}{$\scriptstyle+$}
\relabela <-0.7pt, 1pt> {p3}{$\scriptstyle+$}
\relabela <-0.5pt, 0pt> {1m}{$\scriptscriptstyle 1\!-$}
\relabela <-0.3pt, 0.3pt> {2m}{$\scriptscriptstyle 2\!-$}
\relabela <-0.3pt, 0.3pt> {3m}{$\scriptscriptstyle 3\!-$}
\relabela <0pt, 0pt> {pos1}{$\scriptstyle{\rm positive}$}
\relabela <0pt, 0pt> {pos2}{$\scriptstyle{\rm positive}$}
\relabela <0pt, 0pt> {si1}{$\scriptstyle{\rm side\ of}$}
\relabela <0pt, 0pt> {si2}{$\scriptstyle{\rm side\ of}$}
\relabela <0pt, 0pt> {di1}{$\scriptstyle{\rm disc\ 1}$}
\relabela <0pt, 0pt> {di2}{$\scriptstyle{\rm disc\ 2}$}
\endrelabelbox}
\caption {Example of a positive tiled surface
\label{figure:example of a positive tiled surface}}
\end{figure}

Our work begins with a special case of the embeddability question:
when is a tiled surface which has only positive vertices embeddable in
3--space?  The answer, roughly, is ``most of the time":

For convenience, we call a tiled surface with only positive vertices
a {\em positive} tiled surface.  For such a surface every singularity is
type $aa$ and each singularity is connected to exactly two vertices along
non-adjacent branches of the singular leaves  (see
Figure~\ref{figure:3 tile types}).
Figure~\ref{figure:example of a positive tiled surface}(a) is an example of a
positive tiled surface.  The example is very simple, and so it's easy to
understand the embedding in 3--space which is given in Figure \ref{figure:example
of a positive tiled surface}(b). The surface is
depicted as a Seifert surface for the closed braid $\sigma_1^3$, in Artin's
well-known generators of the braid group. The two discs have been arranged as
concentric discs in 3--space, with disc 2 above disc 1. The two discs are joined
by three half-twisted bands. The singularities in all 3 bands are negative. The
boundary is the negative trefoil knot.

\begin{lemma} 
\label{lemma:embeddability of positive tiled surfaces}
Let $\cF = (F,G,C)$ be a positive tiled surface.  Assume that the
combinatorial data $C$ is subject to a single restriction: the cyclic order
of the singularities around each vertex of valence $\geq 3$ is
counterclockwise when viewed on the positive side of $\cF$. Then $\cF$
is embeddable.
\end{lemma} 

\pf Clearly, in an embedded surface the singular leaves meeting at a
vertex are circularly ordered because the ordering is given precisely
by their order around the braid axis.  Thus the order condition of the
lemma is necessary.  We need to show it suffices when there are no
negative vertices in $G$.
Assuming the order condition is met about each vertex, we shall
construct an embedded orientable surface in three space whose boundary
is a braid, and whose graph of singular leaves in the associated braid
foliation is isomorphic to $G$ with an isomorphic embedding and
isomorphic combinatorial data.

Let $v_1,v_2,\ldots ,v_P$ be the vertices of $\cF$, and let $s_1,s_2,\ldots
s_{S}$ be the singularities, written in order.  Since \cF has no
negative vertices, each singularity is adjacent to exactly two
vertices with the other two singular leaves going to $\partial\cF$.
Let $\delta_i$ be a disc parallel to the $xy$--plane,
centered on the $z$--axis, height $i$, and radius $1/i$.  If $s_i$ is
adjacent to $v_j$ and $v_k$, then connect the discs
$\delta_j$ and $\delta_k$ with small twisted bands at angle $2\pi
i/S$.  The twisted band can twist in either of two ways.  Choose the
twist so that the part connected to the positive part of $\delta_i$ is
positively oriented with respect to $d\theta $ if the sign $s_i$ is
positive, and choose the other twist if $s_i$ is negative.  Note that
the edges of the twisted band can be made arbitrarily close to
straight lines because $1/i$ is a convex function.

Clearly, then, the surface $F^e$ given by the $\delta_i$ and the twisted
bands is a surface with closed braid boundary and associated graph $G$.
The surface
$F^e$ is orientable since a twisted band always connects the discs so
that the upper sides connect to each other.  
All that remains is to check the signs of the
vertices and singularities.  Clearly all the vertices are positive and
the twists were chosen so that the signs of the singularities would
agree.  Thus the braid foliation on $F^e$ and $\cF$ have the same
combinatorial data, and it then follows from Theorem \ref{theorem:unique
embedding} that  $\cF$ is embeddable.  \endpf

We next consider the question of determining a braid word which describes the
boundary of a positive embeddable tiled surface, ie an embeddable tiled
surface which has only positive vertices in its foliation.  Since isotopic
embeddings  of the same tiled surface can have boundaries that differ by a 
conjugation, the answer can only be determined up to conjugation.  There is a
convenient set of generators for the braid group  known as {\em band
generators}.   They are particularly useful in algorithmic questions, having
been to give fast solutions to the word problem in $B_n$ \cite{BKL} and the
conjugacy problem in $B_3$ in
\cite{Xu} and
$B_4$ in \cite{KKL}.

\begin{figure}[htb!]
\centerline{\small
\relabelbox
\epsfxsize 1.7in\epsfbox{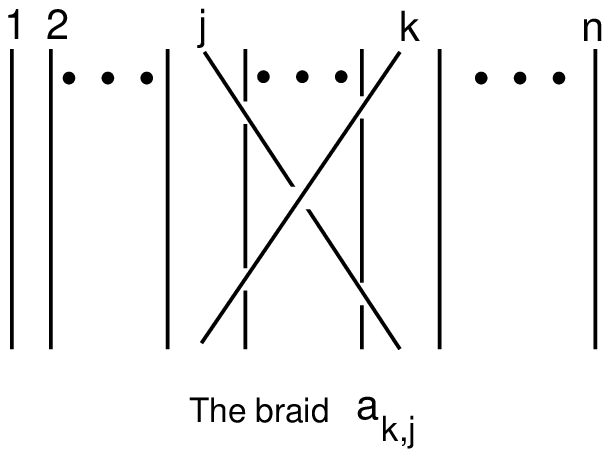}
\extralabel <-44.5mm, 30mm> {1}
\extralabel <-42mm, 30mm> {2}
\extralabel <-30.5mm, 30mm> {$j$}
\extralabel <-17mm, 30mm> {$k$}
\extralabel <-34mm, 0pt> {The braid $a_{k,j}$}
\extralabel <-3mm, 30mm> {$n$}
\endrelabelbox}
\caption{Band generators for the braid group
\label{figure:band generators}}
\end{figure}

Let $k,j$ be integers with $n\geq k\geq j\geq 1$. Let $a_{k,j}$ denote
an elementary braid in which strands $k$ and $j$ exchange places with
strand $k$ crossing over strand $j$ and with both passing in front of
all intermediate braid strands. See Figure
\ref{figure:band generators}. The collection of elementary braids
$a_{k,j}$ and their inverses clearly generate $B_n$ because they
contain as a subset the Artin generators $\sigma_i = a_{i+1,i}$. These
are the {\em band generators} of the braid group.  (See
Appendix~\ref{appendix:testing for conjugacy} for a discussion of
these generators and their relations.)

\begin{lemma}
\label{lemma:braid word for a positive embeddable tiled surface} 
Let \cF be a positive embeddable tiled surface with $P$ positive
vertices (so its boundary is a braid with $P$ strands) and $S$
singularities at angles $\theta_1, \theta_2,\ldots ,\theta_S$ of signs
$\epsilon_1, \epsilon_2, \ldots ,\epsilon_S$.  Notice that every
singularity in the embeddable tiled surface has a exactly two
branches connected to positive vertices. For the singularity which is
at angle $\theta_i$, let $v_{j_i},v_{k_i}$ be the vertices
associated to this singularity, where $k_i>j_i$.  Then the
closed braid given by $B(\cF) =
\prod_{i=1}^{S-1}a^{\epsilon_i}_{k_i,j_i}$ is a representative of
$\partial\cF$. 
\end{lemma}

\pf The proof follows directly from the construction which 
was given in the proof of Lemma~\ref{lemma:embeddability of positive
tiled surfaces}.  That lemma constructed an embedded surface with the
same tiling as \cF.  We abuse notation slightly and call the embedded
surface \cF as well.

 From the construction we know that \cF is made of discs $\delta_i$
centered on the braid axis with twisted bands between these discs.  At
$\theta_i$ there is a band between $\delta_{j_i}$ and $\delta_{k_i}$
with a twist corresponding to the sign $\eps_i$ of the singularity at
$\theta_i$.  The boundary in a neighborhood of $\theta_i$ is then
exactly given by the band generator $a_{k_i,j_i}^{\eps_i}$.
Thus the full closed braid is
$B(\cF) = \prod_{i=1}^{S-1}a^{\epsilon_i}_{k_i,j_i}$.
\endpf

\subsection{The general case: finding the boundary word}
\label{subsection:the boundary word}
We proceed to the general case, where both positive and negative vertices
occur.  The most efficient way to
proceed is to bypass (for the moment) the question of how to test for
embeddability, and assume that we have been given an embeddable tiled surface.

\med{\bf Definition}\qua Let \cF be an embeddable tiled surface with
$P$ positive vertices $v_1,v_2,\ldots ,v_{P}$, in that order on $A$,
and $N$ negative vertices and $S$ singularities at angles $\theta_1,
\theta_2,\ldots ,\theta_S$ of signs $\epsilon_1, \epsilon_2, \ldots
,\epsilon_S$.  Let $\mu$ be the number of components in $\partial\cF$. Recall
that every singularity in the foliation has exactly two
branches connected to positive vertices. For the singularity which is
at angle $\theta_i$, let $v_{j_i},v_{k_i}$ be the orientable surfaces
which are associated to these two vertices, where $k_i>j_i$.  Let
$EB(\cF)$ be the closed braid given by
$$EB(\cF) = \prod_{i=1}^{S-1}a^{\epsilon_i}_{k_i,j_i}.$$
This is a link of one or more components, in 3--space.  The word $EB(\cF)$ is
called the {\em extended boundary word} of $\cF$.  If the tiled surface $\cF$
has only positive vertices, then by Lemma 
\ref{lemma:braid word for a positive embeddable tiled surface}
$EB(\cF) = B(\cF)$.  The word $EB(\cF)$ is given as a word in the band
generators of the braid group
$B_P$, where
$P$ is the number of positive vertices in the tiling. 
Our first lemma tells us that $\partial F$ is represented 
by a word in the braid group $B_{P-N}$.

\begin{lemma}
\label{lemma:n=P-N}
Let $\cF$ be an embeddable tiled surface which has $P$ positive and
$N$ negative vertices. Then the braid index of $\partial F$ is $n = P
- N$.
\end{lemma}
\pf The braid index $n$ is the linking number of $K = \partial F$ with the
braid axis $A$.  Linking number may also be computed as the algebraic
intersection number of $A$ with a surface which $K$ bounds, ie $P-N$. 
\endpf

The next theorem tells us that one of the components of the closed
braid $EB(\cF)$ represents the boundary of the surface $F$.

\begin{theorem} 
\label{theorem:boundary word}
Let $\cF=(F,G,C)$ be an embeddable tiled surface with connected
boundary, $P$ positive vertices and $N$ negative vertices.  Let $K'$
be the link represented by the extended boundary word $EB(\cF)$.  Then
$K'$ is a link with $N+1$ components, at least $N$ of which are
closed 1--braids which are geometrically unlinked from the other
components of $K'$.  Let $K$ be the link which is obtained from $K'$
after deleting $N$ 1--braid components of $K'$. Then $K$ is a
$(P-N)$--braid whose closed braid has 1 component, and this
component represents $\partial F$.
\end{theorem}
\pf Let $\cF'=(F',G',C')$ be the tiled surface induced by removing a
small neighborhood about each negative vertex in $F$ and deleting the
corresponding vertices from $G$ and $C$.  The surface
$F'$ is $F$ {\em with N holes}. The tiled surface $\cF'$ is clearly
embeddable for the following reason: It a subset of the embeddable
tiled surface $\cF$.  Also $\cF'$ has no negative vertices, so Lemma 
\ref {lemma:braid word for a positive embeddable tiled surface} applies,
and $\partial F'$ is represented by
$B(\cF')$.  Thus $K'=\partial\cF'$ is described by the word $B(\cF')$.
Notice that $B(\cF')$ is identical with $EB(\cF)$, but not with
$B(\cF)$.

Thinking of $F'$ as a subset of $F$ embedded in 3--space, we see that
the boundary link of $F'$ contains $N$ small circles which bound discs
each containing a single negative vertex in $F$.  These discs are
disjoint from $F'$ (except on the boundary, of course), thus these $N$
components of $\partial F'$ are geometrically unlinked from the other
components.  Except for these $N$ components, the boundaries of $F$
and $F'$ are identical.  Deleting these $N$ components from $\partial
F'$ yields exactly $\partial F$.
\endpf

\begin{figure}[htb!]
\centerline{\relabelbox\small
\epsfxsize 2.5in\epsfbox{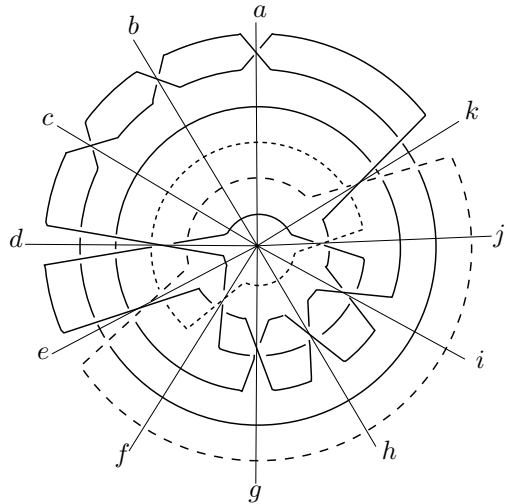}
\relabela <0pt, -1pt> {a}{$a$}
\relabela <0pt, 0pt> {b}{$b$}
\relabela <-1pt, 0pt> {c}{$c$}
\relabela <-4pt, 0pt> {d}{$d$}
\relabela <-2pt, 0pt> {e}{$e$}
\relabela <-2pt, -2pt> {f}{$f$}
\relabela <-1pt, -2pt> {g}{$g$}
\relabela <0pt, -3pt> {h}{$h$}
\relabela <0pt, -2pt> {i}{$i$}
\relabela <0pt, 0pt> {j}{$j$}
\relabela <0pt, 0pt> {k}{$k$}
\endrelabelbox}
\caption{The extended boundary word for the example in Figure
\ref{figure:an example of an embeddable tiled surface}}
\label{figure:boundary word}
\end{figure}

\med{\bf Example}\qua We illustrate Theorem \ref{theorem:boundary word},
using the example in Figure \ref{figure:an example of an embeddable
  tiled surface}. The singularities are at $a,b,c,d,e,f,g,i,j,k$. Of
those, only the singular leaves at $a,b,c,g,h,i$ have two endpoints on
$\partial F$.  Assuming that our tiled surface is embeddable, we
determine its extended boundary word.  The tiling has 6 positive
vertices, 2 negative vertices and 11 singularities. The extended
boundary word $EB(\cF)$ is a 6--braid of length 11 in the band
generators.  It is:
$$
EB(\cF)=a_{2,1}^{3}a_{6,1}^{-1}a_{5,1}^{-1}a_{6,4}^{-1}a_{5,3}^3a_{6,4}a_{5,1}.$$
Figure \ref{figure:boundary word} shows that it has 3 components, two
of which (sketched as dotted curves) are closed 1--braids which are unlinked
from the rest of the braid and from one-another. (It is {\em not}
unique to this example that one has to check carefully to be sure that
they are unlinked from the rest of the braid and from one-another.)
The third component represents the boundary of the surface of Figure
\ref{figure:an example of an embeddable tiled surface}.  Its braid
index, which is the linking number of the axis $A$ with $\partial F$,
may be computed as the number of times the axis pierces \cF, where
intersections are counted algebraically.  Since $P=6, N=2$ this
algebraic intersection number is $6-2=4$, and indeed we see that
(ignoring the two 1--braids) $\partial F$ is a 4--braid.  The only thing
which is not completely obvious at this time is how to instruct a
computer to find a word in the generators of $B_4$ which represents
$\partial F$ from the 3--component 6--braid $EB(\cF)$.  This will be
discussed briefly at the end of Section~\ref{section:enumerating}, and
in more detail in our paper on the implementation of the algorithm.

\subsection{The general case: testing for embeddability}
\label{subsection:testing for embeddability}
We pass to the question of testing the embeddability of an arbitrary
tiled surface $\cF$.  Since the positive tiled surface $\cF'$ is a
subsurface of $\cF$, and since Lemma \ref{lemma:embeddability of
  positive tiled surfaces} gives a complete test for the embeddability
of $\cF'$, it is clear that our general embeddability test must
include the cyclic order test of Lemma \ref{lemma:embeddability of
  positive tiled surfaces} (see (i) of
Theorem~\ref{theorem:embeddability test}, stated below) and a
corresponding condition on the cyclic order around the negative
vertices.  In view of the proof of Theorem \ref{theorem:boundary
  word}, the remaining obstruction to embedding lies in filling in the
disc neighborhoods of the the negative vertices.  The obstruction must
lie in the $b$--arcs, which are not present in a positive tiled
surface.  To describe the obstruction, we need several definitions.

By our hypothesis, the foliation of $\cF$ is radial about each vertex.
This means that around every vertex there is a leaf which meets the
vertex at the angle $\theta$, for every $\theta\in [0,2\pi]$.  Suppose
the singular leaves occur at angles $\theta_1,\dots,\theta_S$.
Consider $b(v_i,v_j)$, a $b$--arc joining vertices $v_i$ and $v_j$.
There is some maximal open interval, $( \theta_m, \theta_n)$ in which
for any $\theta \in ( \theta_m, \theta_n)$, there is a $b$--arc between
vertices $v_i$ and $v_j$ which is homotopic to $b(v_i,v_j)$ rel
endpoints.  Let $[b(v_i,v_j)]$ be the equivalence class given by these
$b$--arcs.  

By a slight abuse of notation we say that the $b$--arc $b(v_i,v_j)$
{\em exists in the $\theta $ interval $(\theta_{p-1},\theta_p)$} if
$(\theta_{p-1},\theta_p) \subset ( \theta_m, \theta_n)$, ie, if some
representative of the equivalence class $[b(v_i,v_j)]$ exists between
angles $\theta_{p-1}$ and $\theta_p$.  Define $gb(v_i,v_j)$ to be a
$gb$--{\em arc} (or {\em generalized} $b$--arc) if either it is a true
$b$--arc, or if $v_i$ and $v_j$ are the positive vertices associated to
a single type $aa$ singularity at $\theta_p$.  In the latter case, we
define the arc $gb(v_i,v_j)$ to lie in then interval
$(\theta_{p-1},\theta_p)$.  Notice that we do {\em not} include the
corresponding arcs for an $ab$ singularity, we will not need them.
When we do not need to distinguish between the $gb$--arcs which are
$b$--arcs and those which are not $b$--arcs, we will use the simpler
notation $v_iv_j$.

\med{\bf Example}\qua Table \ref{table:$gb$-arcs} illustrates the table of
$gb$ arcs for the example of Figure~\ref{figure:an example of an embeddable
  tiled surface}.  The dotted entries indicate the intervals
which end at an $ab$--singularity; for such a singularity there is no
$gb$--arc which is not a $b$--arc.  In a more complicated example the
same would be true for $bb$--singularities.  It is a consequence of our
definitions that there are exactly $N=2$ arcs of type $b$ in every
interval and either one or no arcs which have type $gb$ but not type
$b$. 

\begin{table}
\centerline{\begin{tabular}[hbt!]{l||ccccc}
interval  && $gb$--arcs \\
\hline
$(\theta_1,\theta_2)$ & $b(v_4,v_{2.2})$ & $b(v_5,v_{2.1})$ & $gb(v_1,v_2)$                     
\\
$(\theta_2,\theta_3)$ & $b(v_4,v_{2.2})$ & $b(v_5,v_{2.1})$ &
$gb(v_1,v_2)$                           \\
$(\theta_3,\theta_4)$ & $b(v_4,v_{2.2})$ & $b(v_5,v_{2.1})$ &
$gb(v_1,v_6)$                           \\
$(\theta_4,\theta_5)$ & $b(v_4,v_{2.2})$ & $b(v_5,v_{2.1})$ & $\cdots$  \\
$(\theta_5,\theta_6)$ & $b(v_4,v_{2.2})$ & $b(v_1,v_{2.1})$ & $\cdots$  \\
$(\theta_6,\theta_7)$ & $b(v_6,v_{2.2})$ & $b(v_1,v_{2.1})$ &
$gb(v_5,v_3)$                           \\
$(\theta_7,\theta_8)$ & $b(v_6,v_{2.2})$ & $b(v_1,v_{2.1})$ &
$gb(v_5,v_3)$                           \\
$(\theta_8,\theta_9)$ & $b(v_6,v_{2.2})$ & $b(v_1,v_{2.1})$ &
$gb(v_5,v_3)$                           \\
$(\theta_9,\theta_{10})$ & $b(v_6,v_{2.2})$ & $b(v_1,v_{2.1})$ & $\cdots$  \\
$(\theta_{10},\theta_{11})$ & $b(v_4,v_{2.2})$ & $b(v_1,v_{2.1})$ & $\cdots$ \\
$(\theta_{11},\theta_1)$ & $b(v_4,v_{2.2})$ & $b(v_5,v_{2.1})$ &
$gb(v_1,v_2)$ 
\end{tabular}}
\caption{The table of $gb$--arcs for the example in Figure
\ref{figure:an example of an embeddable tiled surface}
\label{table:$gb$-arcs}}
\end{table} 
Our embeddability test is given by the following theorem:
\begin{theorem} 
\label{theorem:embeddability test} 
Let $\cF$ be a tiled surface whose regions have been labelled in the
manner just described.  Then $\cF$ is embeddable if and only if:
\begin{enumerate}
\item [\rm(i)] The singularities about each positive (respectively  negative)
vertex are positively (respectively  negatively) cyclically ordered in the
fibration, with respect to increasing polar angle $\theta$.
\item [\rm(ii)] The vertices about each positive (respectively  negative)
singularity are positively (respectively  negatively) cyclically ordered on
the oriented braid axis, and
\item [\rm(iii)] The endpoints of a $gb$--arc in the interval $(i-1,i)$
never separate the endpoints of a $b$--arc in the same interval.
\end{enumerate} 
\end{theorem}

\pf  \ We begin the proof by establishing a set of tests which look much more
complicated than the tests in Theorem \ref{theorem:embeddability test}, but
which will turn out to be equivalent to them.   Let $\cF$ be a  tiled surface
with singularities at angles
$\theta_1,\ldots \theta_S$.  We claim that $\cF$ is embeddable if and only if it
passes the following four tests.
\begin{enumerate}
\item \label{cond:ordered sings} The singularities about each positive
  (respectively  negative) vertex are positively (respectively  negatively)
  cyclically ordered in the fibration.

\item \label{cond:ordered verts} The vertices about each positive
  (respectively  negative) singularity are positively (respectively  negatively)
  cyclically ordered in the braid axis.  

\item \label{cond:b-arc separation} Let $vw$ be a $b$--arc which exists
  during the $\theta$--interval $(i-1,i)$.  Then $\cF$ is not
  embeddable if the vertices $v$ and $w$ are separated on $A$ by the
  endpoints of any other $b$--arc which exists in the $\theta$--interval
  $(i-1,i)$.
  
\item \label{cond:aa separation} Suppose the singularity at
  $\theta_{i}$ is type $aa$, between $a$--arcs at vertices $v$ and $w$.
  Then $\cF$ is not embeddable if the vertices $v$ and $w$ are
  separated on $A$ by the endpoints of one of the $b$--arcs which exist
  in the $\theta$--interval $(i-1,i)$.
  
\item \label{cond:ab separation} Suppose the singularity at
  $\theta_{i}$ is type $ab$, between an $a$--arc with vertex endpoint
  $x$ and a $b$--arc $uv$, where $u$ is positive.  Then $\cF$ is not
  embeddable if there is different $b$--arc, say $yz$, which occurs
  during $(i-1,i)$ such that $x$ is separated from $uv$ on $A$ by
  $yz$.  

\item \label{cond:bb separation} If the singularity at $\theta_{i}$ is
  type $bb$, let $u,v,w,x$ be the vertices of the $bb$ tile $T$,
  oriented as they are encountered in traversing $\partial T$
  counterclockwise on the positive side of $\cF$, with $u$ positive.
  Then $\cF$ is not embeddable if there is a $b$--arc in the interval
  $(i-1,i)$ which separates $uv$ from $wx$.
\end{enumerate}

To prove the claim we first notice that (\ref{cond:ordered sings}) is
a necessary condition for the surface to be embeddable, because the
foliation is radial in a sufficiently small neighborhood of every
vertex.  Similarly, (\ref{cond:ordered verts}) is necessary because
the singular leaves meeting at a singularity are embedded in a single
$H_\theta$.  The vertices at the ends of the leaves are all in the
braid axis on the boundary of $H_\theta$ inducing an order on the
vertices.  The sign of the singularity indicates whether the surface
is locally oriented compatibly with $H_\theta$, or with reversed
orientation. 

Consider the intersections of the given orientable surface $F$ with
the fibers of $\cH$, as $\theta$ ranges over the interval $[0,2\pi]$.
Let $N$ be the number of negative vertices and let $P$ be the number
of positive vertices.  An Euler-characteristic count shows that there
must be $P$ arcs in all, with exactly $N$ of them type $b$ and the
remaining ones type $a$.

We first find necessary conditions for embeddability.  If $F$ is
embedded, then it has no self-intersections.  Since $F$ intersects
each non-singular fiber transversally, it follows that a necessary
condition for embeddability is that $F\cap H_\theta$, where
$H_\theta$ is non-singular, be a collection of pairwise disjoint arcs,
with $N$ of them of type $b$ and the remaining $P-N$ type $a$.  The
$b$ arcs divide $H_\theta$, but the $a$ arcs do not.  If there are
$b$--arcs, say $uv,wx\subset H_\theta$, they will intersect if and only
if $u, v$ separate $w,x$ on $A=\partial H_\theta$.  Let
$H_{\theta_1},\dots,H_{\theta_s}$ be the singular fibers, in their
natural cyclic order in the cycle of fibers around $A$, with
subscripts mod $s$.  If $\theta$ varies over the open interval
$\theta\in (\theta_{i-1},\theta_i)$ its intersections with $F$ will
be modified by isotopy rel $A$.  Thus a necessary condition for $\cF$
to be embeddable is that it pass test (\ref{cond:aa separation}) for
some $\theta\in (\theta_{i-1},\theta_i)$ for every $i=1,2,\dots,S$.

Next we ask what happens to the intersections of our embedded
orientable surface $F$ with $H_\theta$ when $\theta$ passes through
a singular angle in the fibration.  There are three types of
singularities, ie  $aa$, $ab$, and $bb$.  It's easy to see that the
arcs in the set $F\cap H_\theta$ only change in a manner which can
be realized by an isotopy after an $aa$ singularity, but that is not
the case after an $ab$ or $bb$ singularity.  As we approach the
singular fiber which separates the intervals $(i-1,i)$ and $(i,i+1)$
during an $aa$ singularity the $a$--arcs with endpoints at $v,w$ must
approach one-another.  But if $v$ and $w$ are separated on $A$ by the
endpoints of a $b$ arc which exists during the interval $(i-1,i)$ that
will be impossible without a self-intersection in $F$.  The reasons
are the same for type $ab$ and $bb$.  Thus, tests (\ref{cond:aa
  separation})--(\ref{cond:bb separation}) are also necessary
conditions for embeddability.

In fact these tests are also sufficient.  Assume that all 6 tests have
been passed.  Then $F\cap H_\theta$ is a collection of pairwise
disjoint arcs, with $N$ of them type $b$ and $P-N$ of them type $a$,
for every non-singular fiber.  Also, in every singular fiber there is
exactly one pair of intersecting arcs, namely the leaves of the
associated saddle-point tangency.  The union of all of the arcs
$F\cap H_\theta$ as $\theta$ varies over the closed interval
$[0,2\pi]$ is the trace of the isotopy of $F\cap H_\theta$ as
$\theta$ varies over $[0,2\pi ]$.  The claim is proved.   

To complete the proof of Theorem \ref{theorem:embeddability test} we
now observe that test (i) of the theorem is identical with
(\ref{cond:ordered sings}) of this lemma, and test (ii) is identical
with (\ref{cond:ordered verts}).  Test (iii) of the theorem is
identical with (\ref{cond:b-arc separation}) plus (\ref{cond:aa
  separation}).  It remains to show that  tests
(\ref{cond:ab separation}) and (\ref{cond:bb separation}) are subsumed
by test (iii).

The key fact to notice is that the changes as we pass through an
$ab$ (respectively   $bb$) singularity at the angle $\theta_i$ involve exactly
one (respectively  two) $b$ arc (respectively  arcs).  All other $b$--arcs in the
interval $(i,i+1)$ are identical with those in the interval $(i-1,i)$.

Consider test (\ref{cond:bb separation}) first.  Suppose that there is
a $bb$ singularity at $\theta_i$, as in (\ref{cond:bb separation}),
with $b$ arcs $uv$ and $wx$ in $(i-1,i)$, and new $b$--arcs $vw$ and
$ux$ in $(i,i+1)$.  Suppose also that there is a $b$--arc $cd$ in the
interval $(i-1,i)$ which separates $uv$ and $wx$.  Then $c$ and $d$
separate $u$ and $v$.  However, in the interval $(i,i+1)$ (ie  after
the singularity) there will be new $b$ arcs $vw, ux$.  The $b$--arc
$cd$ will still be present.  But that is impossible by
(\ref{cond:b-arc separation}), because $c$ and $d$ separate $u$ and
$v$.

Consider test (\ref{cond:ab separation}) next.  Suppose that there is
an $ab$ singularity at $\theta_i$, as in (\ref{cond:ab separation}),
such that the $b$--arc $uv$ is in $(i-1,i)$ and the $b$--arc $xv$ is in
$(i,i+i)$.  All other $b$--arcs in $(i-1,i)$ are also in $(i,i+1)$.
Suppose that $x$ is separated from the $b$--arc $uv$ by some other
$b$--arc $yz$ in the interval $(i-1,i)$.  By
(\ref{cond:b-arc separation}),
the arc $yz$ cannot cross $uv$.  This means that $y$ and $z$ separate
$x$ from both $u$ and $v$.  Passing to the interval $(i,i+1)$, the arc
$yz$ is still present.  However $yz$ crosses $xv$.  But that is also
impossible, by 
(\ref{cond:b-arc separation}).  The proof of Theorem
\ref{theorem:embeddability test} is complete.  \endpf

\med{\bf Example}\qua We illustrate the embeddability test on the example
which was given in Figure \ref{figure:an example of an embeddable tiled
  surface}, using the data in Table \ref{table:$gb$-arcs}.   Recall
that the direction of increasing polar angle $\theta$ is
counterclockwise (respectively  clockwise) about a positive (respectively  negative)
vertex.  An easy check shows that the order is correct about every
vertex, so the example passes test (i) of Theorem
\ref{theorem:embeddability test}.  Similarly, test (ii) is passed.  We
turn to test (iii).   There are two negative vertices, and so there are
two $b$--arcs in each non-singular fiber.  Inspecting Figure
\ref{figure:an example of an embeddable tiled surface}, we
see that there is a $b$--arc joining vertices $v_4$ and $v_{2.2}$ in
the interval $(\theta_{10},\theta_6)$, and one joining vertices $v_6$
and $v_{2.2}$ in the interval $(\theta_6,\theta_{10})$.   These are the
entries in the first column of Table \ref{table:$gb$-arcs}.   Similarly
for the $b$--arcs which end at vertex $v_{2.2}$, which are recorded in
the second column.   As for the $gb$--arcs, we see that there are
$aa$--singularities at
$\theta_1,\theta_2,\theta_3,\theta_4,\theta_7,\theta_8, \theta_9$,
explaining the entries in the third column of Table
\ref{table:$gb$-arcs}.  Inspecting the rows of the table, one at a
time, we verify that the endpoints of a $gb$--arc (remembering that
$gb$--arcs include $b$--arcs) never separate the endpoints of a $b$--arc.
Thus test (iii) of Theorem \ref{theorem:embeddability test} is also
passed. 

\subsection{Eliminating inessential tiled surfaces}
\label{subsection:eliminating inessential surfaces}
For efficiency, we will want to be able to eliminate inessential tiled
surfaces as we do the enumeration.  The following proposition will allow us to do
so.
\begin{prop}
\label{proposition:inessential b-arcs}
Let $\cF=(F,G,C)$ be a tiled surface which passes the embeddability
test of Theorem \ref{theorem:embeddability test}.  Then $\cF$ is
essential if, for every $b$--arc $\beta$ in the tiled surface, the two
points in $\partial\beta$ are not adjacent in the cyclic ordering of
the vertices on $A$.
\end{prop}

\pf  
Let $\beta$ be a $b$--arc in $F\cap H_\theta$ for some
non-singular fiber $H_\theta$.  Recall that $\beta$ divides divides 
$H_\theta$ into two subdiscs, $\delta_1$ and
$\delta_2$, and that $\beta$ is {\em inessential} if one of these
discs, say $\delta_1$, has empty intersection with $K$.  The subdisc
$\delta_1$ may contain other $b$--arc within it, but we may assume that
$\beta$ is an innermost $b$--arc and $\delta_1$ contains no other $b$--arcs.
  
In this case, we can simplify the braid foliation of $\cF$ by pushing $F$
through a neighborhood of $\delta_1$ in 3--space to eliminate two
points of $A\cap F$.  We can detect this situation combinatorially
because, if $\beta$ is inessential and innermost, its two endpoints
will be adjacent vertices on $A$.  \endpf

\section{Enumerating closed braid representatives of the unknot}
\label{section:enumerating}
Our task in this section is develop a procedure for enumerating the
closed $n$--braid representatives of the unknot, up to conjugacy, for
each fixed $n$.  From now on we restrict our attention to the
case when the $F$ is a disc. To stress this, we use the notation $\cD
= (D,G,C)$.
Each closed $n$--braid representative of the unknot is the boundary of
a embeddable tiled disc, and our plan is to enumerate all embeddable
tiled discs and read off their boundary words.  If $n\geq 4$ there
will be infinitely many conjugacy classes.  Our measure of complexity
for the enumeration is the number of vertices in the embeddable tiled
disc, ie  the number of points in $A\cap D$.

\begin{figure}[htb!]
\centerline{\relabelbox\small
\epsfxsize 3.5in\epsfbox{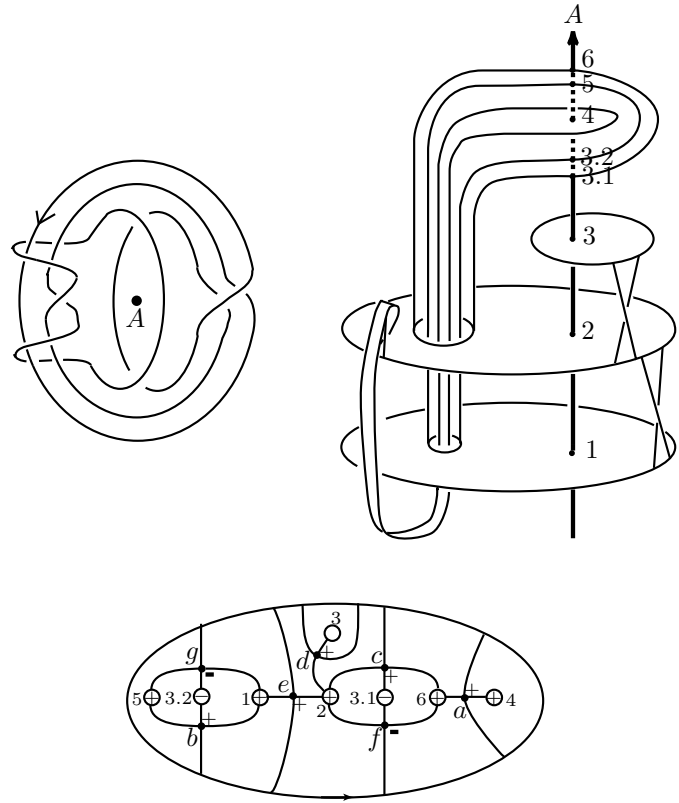}
\relabela <-0.5pt, 0pt> {a}{$a$}
\relabela <0pt, -0.5pt> {b}{$b$}
\relabela <-0.5pt, 0pt> {c}{$c$}
\relabela <-1pt, -1pt> {d}{$d$}
\relabela <0pt, 0pt> {e}{$e$}
\relabela <-3pt, 0pt> {f}{$f$}
\relabela <0pt, 0pt> {g}{$g$}
\relabela <-2pt, 0pt> {A1}{$A$}
\relabela <-2pt, 0pt> {A2}{$A$}
\relabela <0pt, 0pt> {1}{$1$}
\relabela <0pt, 0pt> {2}{$2$}
\relabela <0pt, 0pt> {3}{$3$}
\relabela <0pt, 0pt> {311}{$3.1$}
\relabela <-1pt, 0pt> {312}{$\scriptstyle3.1$}
\relabela <0pt, 0pt> {321}{$3.2$}
\relabela <-1pt, 0pt> {322}{$\scriptstyle3.2$}
\relabela <0pt, 0pt> {4}{$4$}
\relabela <0pt, 0pt> {5}{$5$}
\relabela <0pt, 0pt> {6}{$6$}
\relabela <-0.5pt, 0pt> {11}{$\scriptstyle1$}
\relabela <0pt, 0pt> {21}{$\scriptstyle2$}
\relabela <0pt, 0pt> {31}{$\scriptstyle3$}
\relabela <-0.5pt, 0pt> {41}{$\scriptstyle4$}
\relabela <-1pt, 0pt> {51}{$\scriptstyle5$}
\relabela <-0.5pt, 0pt> {61}{$\scriptstyle6$}
\relabela <-2.6pt, .7pt> {m1}{$\scriptstyle-$}
\relabela <-2.6pt, .7pt> {m2}{$\scriptstyle-$}
\relabela <-2.3pt, -.5pt> {p1}{$\scriptstyle+$}
\relabela <-2.1pt, -.3pt> {p2}{$\scriptstyle+$}
\relabela <-2.1pt, -.3pt> {p3}{$\scriptstyle+$}
\relabela <-2pt, -.5pt> {p4}{$\scriptstyle+$}
\relabela <-1.5pt, .5pt> {p5}{$\scriptstyle+$}
\relabela <-2.1pt, -.3pt> {p6}{$\scriptstyle+$}
\relabela <-2.1pt, -.7pt> {p7}{$\scriptstyle+$}
\relabela <-1.8pt, -.5pt> {p8}{$\scriptstyle+$}
\relabela <-2.3pt, -.5pt> {p9}{$\scriptstyle+$}
\relabela <-2.3pt, -.5pt> {p10}{$\scriptstyle+$}
\endrelabelbox}
\caption{Several views of a closed braid which bounds an interesting disc
\label{figure:interesting disc}}
\end{figure}

We refer the reader to Figure \ref{figure:interesting disc}.   
The top left sketch is a closed 4--braid diagram for the unknot.  It is
readily seen to be the unknot, however the disc that it bounds is a
little obscure.  The top right sketch illustrates the same closed
braid from a different angle.  The disc meets the braid axis $A$ in 8
points (the vertices of the induced foliation), labelled
$1,2,3,3.1,3.2,4,5,6$.  The two labelled 3.1 and 3.2 are negative
vertices.  The singularities are not labelled, but there are 7 of them.
(By an Euler characteristic argument, there is always exactly one
fewer singularity than vertex.)  Singularities occur on each of the
three narrow twisted bands, and each of the narrow, vertical tubes
coming up out of discs 1 and 2 each have two singularities of opposite
signs.  Vertices are labelled with numbers, singularities with letters.
Both singularities and vertices have a label of plus or minus, too.
There is an embedding of a model disc \cD in 3--space which realizes
this geometry.

The bottom sketch in Figure \ref{figure:interesting disc} shows the
model disc.  The graph of the singular leaves has been pulled back to
the model disc, and decorated to show the order and signs of the
vertices and singularities.  This bottom sketch is our embeddable
tiled disc.  The problem which we address now is this: let us suppose
that we were handed the 4--braid example $K$ which is illustrated in
the left sketch in Figure \ref{figure:interesting disc}, and we want
to verify algorithmically that it is the unknot.  Our plan is to
enumerate a suitably long list of foliated discs whose boundaries are
4--braids, and to check our given example $K$ against the members of
the list.  So we need to learn how to generate, systematically, a list
of embeddable tiled discs all of whose boundaries have braid index 4,
which is long enough to contain the example in Figure
\ref{figure:interesting disc}.

In view of Lemma \ref{lemma:n=P-N} our plan
is to fix $n=P-N$ and to enumerate embeddable tiled discs in order of
increasing $v$.  This is the same as enumerating embeddable tiled
discs with $(P,N) = (n,0),(n+1,1),(n+2,2),\dots$ in that order.

To enumerate the embeddable tiled discs we apply ideas first used in
the proof of Lemma 1 of \cite{Bennequin}.  Again we refer the reader
to \cite{BF} for details, and give a brief summary here.  To do so we
introduce a move which is guided by the foliation of the surface and
allows us to change an arbitrary embeddable tiled disc to a new
embeddable tiled disc.  This new embeddable tiled disc is simpler in
the sense that it has fewer negative vertices.  Our algorithm will
attempt to reverse this process, starting with a simple embeddable
tiled disc and generating more complex ones.

\begin{figure}[htb!]
\centerline{\relabelbox\small
\epsfxsize 4in\epsfbox{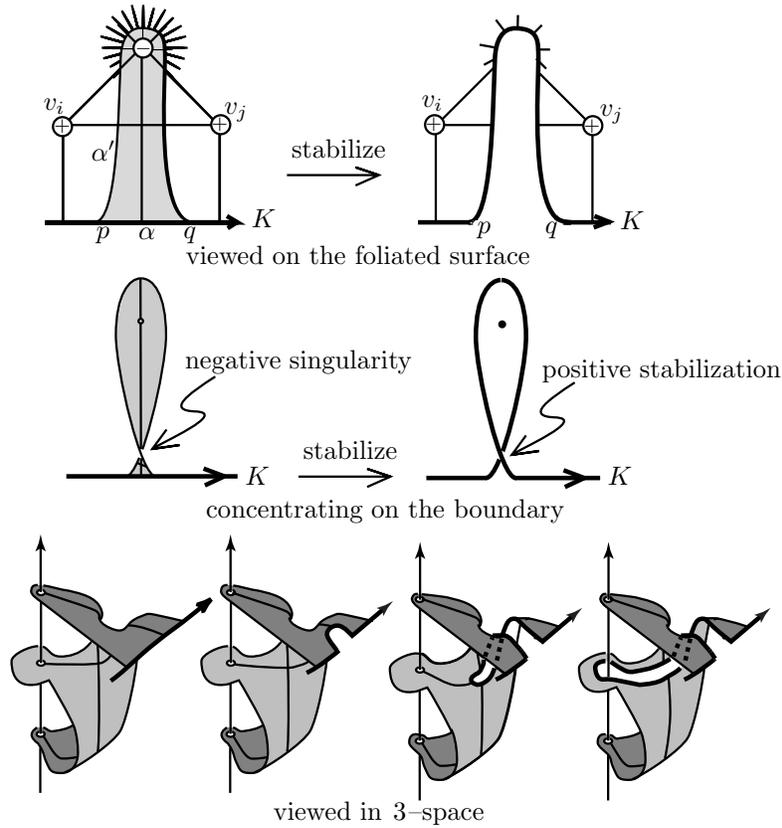}
\relabel {stab1}{stabilize}
\relabel {stab2}{stabilize}
\relabela <-15pt, -4pt> {view1}{viewed on the foliated surface}
\relabel {neg}{negative singularity}
\relabela <0pt, 2pt>  {pos}{positive stabilization}
\relabela <-15pt, 0pt> {conc}{concentrating on the boundary}
\relabel {K1}{$K$}
\relabel {K2}{$K$}
\relabel {K3}{$K$}
\relabel {K4}{$K$}
\relabel {v1}{$v_i$}
\relabel {v2}{$v_j$}
\relabela <-0pt, 2pt> {v3}{$v_i$}
\relabel {v4}{$v_j$}
\relabela <-0pt, -1pt> {a1}{$\alpha$}
\relabela <-4pt, 0pt> {a2}{$\alpha'$}
\relabela <-0pt, -1pt> {pl1}{$p$}
\relabela <-0pt, -1pt> {q1}{$q$}
\relabel {pl2}{$p$}
\relabela <-7pt, 0pt> {q2}{$q$}
\relabela <-1.7pt, 1pt> {p1}{$\scriptstyle+$}
\relabela <-1.4pt, 1pt> {p2}{$\scriptstyle+$}
\relabela <-1.4pt, 1pt> {p3}{$\scriptstyle+$}
\relabela <-1.4pt, 1pt> {p4}{$\scriptstyle+$}
\relabela <-2pt, 1.3pt> {m1}{$\scriptstyle-$}
\relabela <-15pt, 0pt>  {view2}{viewed in $3$--space}
\endrelabelbox}
\caption{Stabilization along an $ab$--tile, viewed from three perspectives
\label{figure:stabilization along an ab-tile}}
\end{figure}

\med{\bf Definition}\qua {\sl Stabilization along an $ab$--singularity}\qua 
Recall that, by hypothesis, the foliation is radial is a neighborhood
of each vertex in the tiling.  The top row in Figure
\ref{figure:stabilization along an ab-tile} shows how, any time there
is an $ab$--singularity in the foliation, we may push $K$ across the
singularity and its associated negative vertex, in a neighborhood of
the separating leaf which meets $K$, to a new position which is again
everywhere transverse to the foliation.  It follows that after we do
this move we will have a new closed braid representative, say
$K^\star$, of the unknot.  Notice that after stabilizing, a $bb$
singularity may have become an $ab$ singularity.   The middle row of
pictures shows why the move {\it increases} the braid index from $n$
to $n+1$, while {\it decreasing} the number $N$ of negative vertices
from $N$ to $N-1$.   The bottom row shows our stabilization move on the
embedded surface in 3--space.   If one looks carefully one can see the
half-twist which has been introduced in the course of the push.  We
note that the pictures of $ab$--tiles in the bottom row of Figure
\ref{figure:stabilization along an ab-tile} are deformations of the
picture in Figure \ref{figure:3 tile types}: we stretched out the top
sheet to make visible a neighborhood of the singular leaf.

\begin{theorem}
\label{theorem:constructing all embeddable tiled discs}
Let $\cD$ be an arbitrary embeddable tiled disc.  Suppose that the
graph of $\cD$ contains $P$ positive vertices and $N$ negative
vertices.  Then there exists a sequence of embeddable tiled discs:
$$\cD = \cD_0\to\cD_1\to\dots\to\cD_{N,}$$
where each $\cD_{i+1}$ is
obtained from $\cD_i$ by a single $ab$--stabilization, 
%%% mdh changed "such that" to "so"
so
$\cD_N$ has only $aa$--singularities.  If 
%%% mdh changed "we start" to "the initial"
the initial
the  \cD is essential, then so is each $\cD_i$.

Equivalently, every embeddable tiled disc with $P$ positive
vertices and $N$ negative vertices may be constructed by finding a
embeddable tiled disc the graph of which is a tree of $P$
positive vertices, then adding $N$ $ab$--tiles, one at a time, to the
graph.  At each addition stage, the new vertex and singularity are
inserted into the orders of the older vertices and singularities and
in such a way that the new graph corresponds to an embeddable tiled
disc.  If the disc to be constructed is essential, then each
intermediate disc may also be assumed to be essential.
\end{theorem}

\pf We begin with the given embeddable tiled disc $\cD$, which, by
hypothesis, contains $N$ negative vertices.  If $N=0$ we are done, so
assume that $N>0$.  From Figure \ref{figure:3 tile types} we can see
that the foliation of $\cD_0$ necessarily contains singularities of
type $bb$ or $ab$, because singularities of type $aa$ only connect to
positive vertices.  But if there are singularities of type $bb$, then
there must also be singularities of type $ab$ because a $bb$ tile can
only be glued to another $bb$ tile or to an $ab$ tile.  However, a
subsurface of $D$ cannot be composed entirely of $bb$--tiles, for if
it were it would be closed, and also entirely in the interior of
$D$, which is absurd.  So we may assume that there is at least one
$ab$--singularity.  It is then possible to stabilize along the $ab$
singularities, one at a time, as in Figure \ref{figure:stabilization
  along an ab-tile} until we obtain a tiled disc $\cD_N$ which has no
negative vertices.

We must
show that each $\cD_i$ is embeddable and has no inessential
$b$--arcs.   Since the graph of $\cD_N$ has no negative vertices, its
singularities must all be type $aa$.   Since \cD is embeddable, and
$\cD_i\subset \cD$ for all $i$, it follows
that each $\cD_i$ is an embeddable tiled disc.   The $b$--arcs of each $\cD_i$
are also $b$--arcs of $\cD$.   Consider
a $b$--arc of $\cD_i$.   It is isotopic to one inside a non-singular
$H_\theta $ which we may assume is sufficiently far from any singular
$H_\phi $.   Since $ab$--stabilization occurs in a neighborhood of a
singularity which is inside a singular $H_\phi $,
$$\partial D \cap H_\theta \subset \partial D_i \cap H_\theta,  $$
and the $b$--arc is essential in $D$ if and only if it is essential in
$D_i$. 

Thus, the sequence of tiled discs is a sequence of (essential)
embeddable tiled discs.  \endpf

In the previous theorem we examined how arbitrary essential tiled
discs can be simplified by stabilization along $ab$--singularities. We showed 
that after some number of such stabilizations one arrives at a positive tiled
disc. The latter contains only positive vertices, and so has only
$aa$--singularities. Its graph of singular leaves is a tree. 
As was shown by Bennequin in \cite{Bennequin}, the new tiled disc
can then be further simplified to the trivially tiled disc by removing
$aa$--singularities which belong to vertices of valence 1, one at a time. 
That process is known as {\it destabilization along
$aa$--singularities}. 
See Figure \ref{figure:destabilization along
aa-singularities}.

\begin{figure}[htb!]
\centerline{\epsfxsize 3in \epsfbox{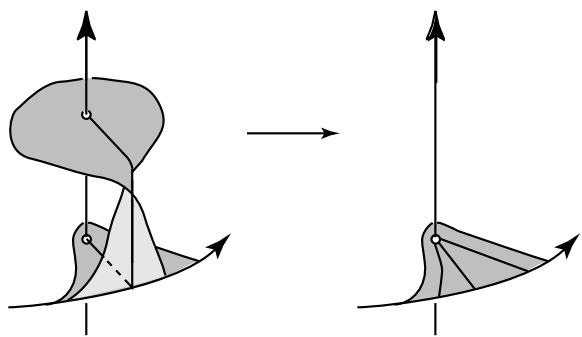}}
\caption{Destabilization along an $aa$--tile having a vertex of valence 1
\label{figure:destabilization along aa-singularities}}
\end{figure} 
  
See \cite{BF} for a full discussion of how the two
processes, ie a finite number of stabilization along $ab$--singularities
followed by a finite number of destabilizations along $aa$--singularities
can be used
to change an arbitrary closed
$n$--braid representative of the unknot to an $(n+m)$--braid representative which
bounds a simpler tiled disc and
thence to the trivial 1--braid representative.  We now consider the reverse of
stabilization along
$ab$--singularities.
%%% mdh
%%%  and destabilization along $aa$--singularities. 
This reverse
proceess will be used to build up all tiled discs. 

\begin{theorem}
\label{theorem:enumerating embeddable tiled discs}
All possible embeddable tiled discs of fixed braid index $n$ may be
enumerated in order of increasing $v$ by the following (not
necessarily efficient) procedure:
\begin{itemize}
\item Enumerate all positive tiled discs with $n$ vertices, testing
  each for embeddability and discarding non-embeddable examples.  This
  gives a list of all possible embeddable tiled discs of braid index
  $n$ with $v=n$ vertices.
\item To enumerate all embeddable tiled discs with $n+2$ vertices,
  first enumerate all positive tiled discs with $n+1$ vertices,
  testing each for embeddability and discarding non-embeddable
  examples.  Then add one $ab$ tile, using the reverse of stabilization along
an $ab$--singularity, in all
  possible ways, testing every position for embeddability and
  essentiality and discarding non-embeddable or inessential examples.
  This produces a list of all possible essential embeddable tiled
  discs of braid index $n$ with $v=n+2$ vertices.
\item To enumerate all embeddable tiled discs with $n+4$ vertices,
  first enumerate all positive tiled discs with $n+2$ vertices, then
  discard any non-embeddable examples.  Then add two
  $ab$ tiles in all possible ways, and
  discard any non-embeddable or inessential
  examples.  This produces a list of all possible embeddable tiled
  discs of braid index $n$ with $v=n+4$ vertices.
\item Continue in this way for $v=n+6, n+8,\dots$.  For $v=n+2i$, the
  enumeration begins with all possible positive embeddable tiled
  discs with $n+i$ vertices, and continues by adding $i$ tiles of type $ab$.
\end{itemize}
\end{theorem}

\pf By Lemma \ref{lemma:n=P-N} we know that $n=P-N$. Since $v=P+N$, the
enumeration begins with $(P,N)=(n,0)$ and continues with
$(n+1,1),(n+2,2),\dots$.  By Theorem
\ref{theorem:constructing all embeddable tiled discs} every essential
embeddable tiled disc will eventually appear on the list.  A key fact
in the enumeration is that 
%%%mdh 
an
embeddable tiled disc can never be
obtained from a non-embeddable tiled disc by the reverse of stabilization along
an $ab$--singularity.  The reason is that each $\cD_i$ in any sequence
of tiled discs constructed by repeated $ab$--stabilization on a
embeddable tiled disc is embeddable.  Similarly, by Lemma
\ref{theorem:constructing all embeddable tiled discs} an embeddable tiled disc
with essential
$b$--arcs can never be obtained from a embeddable tiled disc having
%%% mdh changed "only" "any"
any inessential $b$--arcs by  
the reverse of stabilization along an $ab$--singularity. 
\endpf

\begin{cor}
\label{corollary:enumerating conjugacy classes}
All possible conjugacy classes of $n$--braid representatives of the
unknot may be enumerated in order of complexity $(n,v)$ by enumerating
all embeddable tiled discs, using Theorem \ref{theorem:enumerating embeddable tiled discs}, and then
using Theorem \ref{theorem:boundary word} to determine the associated
boundary words.\endpf
\end{cor}

\begin{figure}[htb!]
\let\ss\scriptstyle
\centerline{\relabelbox\small
\epsfxsize 5in\epsfbox{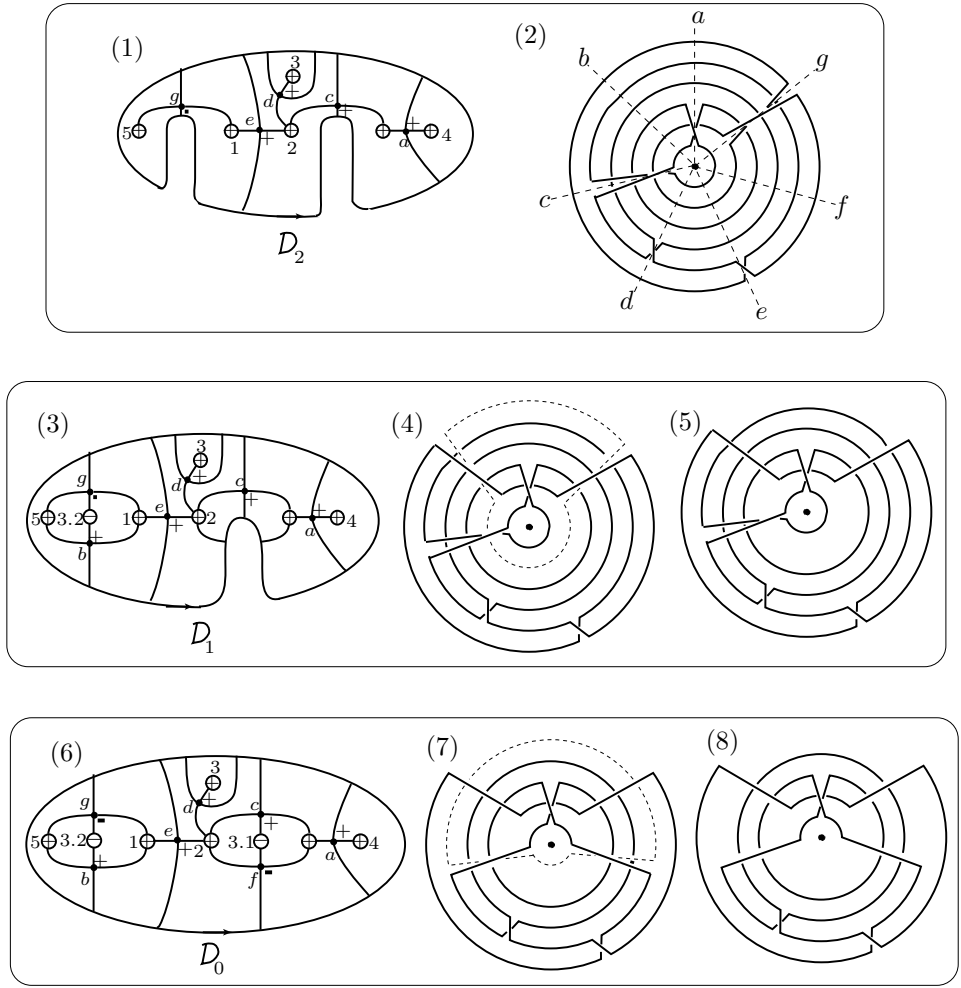}
\relabel {01}{(1)}
\relabel {02}{(2)}
\relabel {03}{(3)}
\relabel {04}{(4)}
\relabel {05}{(5)}
\relabel {06}{(6)}
\relabel {07}{(7)}
\relabel {08}{(8)}
\relabela <-2.5pt, -.3pt> {p1}{$\ss+$}
\relabela <-2.2pt, -.3pt> {p2}{$\ss+$}
\relabela <-2.5pt, -.3pt> {p3}{$\ss+$}
\relabela <-2pt, -.3pt> {p4}{$\ss+$}
\relabela <-.5pt, 1.7pt> {p5}{$\ss+$}
\relabela <-2.5pt, -.3pt> {p6}{$\ss+$}
\relabela <-1.5pt, -1.3pt> {p7}{$\ss+$}
\relabela <-1.5pt, -.3pt> {p8}{$\ss+$}
\relabela <-.5pt, -.8pt> {p9}{$\ss+$}
\relabela <-2.5pt, -.3pt> {p10}{$\ss+$}
\relabela <-2.5pt, -.3pt> {p11}{$\ss+$}
\relabela <-2.5pt, -.3pt> {p12}{$\ss+$}
\relabela <-2.5pt, -.3pt> {p13}{$\ss+$}
\relabela <-2.5pt, -.3pt> {p14}{$\ss+$}
\relabela <-2pt, -.3pt> {p15}{$\ss+$}
\relabela <-1.5pt, .7pt> {p16}{$\ss+$}
\relabela <-2.5pt, -.3pt> {p17}{$\ss+$}
\relabela <-1pt, -1.3pt> {p18}{$\ss+$}
\relabela <-.5pt, -.3pt> {p19}{$\ss+$}
\relabela <-2pt, -.8pt> {p20}{$\ss+$}
\relabela <-2.5pt, -.3pt> {p21}{$\ss+$}
\relabela <-2pt, -.3pt> {p22}{$\ss+$}
\relabela <-2.5pt, -.3pt> {p23}{$\ss+$}
\relabela <-2.5pt, -.3pt> {p24}{$\ss+$}
\relabela <-2.5pt, -.3pt> {p25}{$\ss+$}
\relabela <-2.5pt, -.3pt> {p26}{$\ss+$}
\relabela <-2pt, -.3pt> {p27}{$\ss+$}
\relabela <-1.5pt, .7pt> {p28}{$\ss+$}
\relabela <-1.5pt, -1.3pt> {p29}{$\ss+$}
\relabela <-.5pt, -.3pt> {p30}{$\ss+$}
\relabela <-1.5pt, .2pt> {p31}{$\ss+$}
\relabela <-2.5pt, -.3pt> {p32}{$\ss+$}
\relabela <-2.7pt, 0.8pt> {m1}{$\ss-$}
\relabela <-2.7pt, 0.8pt> {m2}{$\ss-$}
\relabela <-2.7pt, 0.8pt> {m3}{$\ss-$}
\relabela <0pt, 0pt> {a1}{$\ss a$}
\relabela <0pt, 0pt> {a2}{$\ss a$}
\relabela <0pt, 0pt> {a3}{$\ss a$}
\relabela <0pt, -0.5pt> {a4}{$a$}
\relabela <0pt, 0pt> {b1}{$\ss b$}
\relabela <0pt, 0pt> {b2}{$\ss b$}
\relabela <0pt, 0pt> {b4}{$b$}
\relabela <0pt, 0pt> {c1}{$\ss c$}
\relabela <0pt, 0pt> {c2}{$\ss c$}
\relabela <0pt, 0pt> {c3}{$\ss c$}
\relabela <0pt, 0pt> {c4}{$c$}
\relabela <0pt, 0pt> {d1}{$\ss d$}
\relabela <-.5pt, -.5pt> {d2}{$\ss d$}
\relabela <0pt, 0pt> {d3}{$\ss d$}
\relabela <0pt, 0pt> {d4}{$d$}
\relabela <0pt, 0pt> {e1}{$\ss e$}
\relabela <0pt, 0pt> {e2}{$\ss e$}
\relabela <0pt, 0pt> {e3}{$\ss e$}
\relabela <0pt, 0pt> {e4}{$e$}
\relabela <-2pt, 0pt> {f1}{$\ss f$}
\relabela <-2pt, 0pt> {f4}{$f$}
\relabela <0pt, 0pt> {g1}{$\ss g$}
\relabela <0pt, 0pt> {g2}{$\ss g$}
\relabela <0pt, 0pt> {g3}{$\ss g$}
\relabela <0pt, 0pt> {g4}{$g$}
\relabela <-1pt, 1.5pt> {00}{$\ss 0$}
\relabela <-1pt, 0pt> {11}{$\ss 1$}
\relabela <-1pt, 0pt> {12}{$\ss 1$}
\relabela <-1pt, 0pt> {13}{$\ss 1$}
\relabela <-1pt, 1.5pt> {14}{$\ss 1$}
\relabela <-1pt, 0pt> {21}{$\ss 2$}
\relabela <-1pt, 0pt> {22}{$\ss 2$}
\relabela <-1pt, 0pt> {23}{$\ss 2$}
\relabela <-1pt, 1.5pt> {24}{$\ss 2$}
\relabela <-1pt, 0pt> {31}{$\ss 3$}
\relabela <-1pt, 0pt> {311}{$\ss 3.1$}
\relabela <-1pt, 0pt> {32}{$\ss 3$}
\relabela <-1pt, 0pt> {321}{$\ss 3.2$}
\relabela <-2pt, 0pt> {322}{$\ss 3.2$}
\relabela <-1pt, 0pt> {33}{$\ss 3$}
\relabela <-1pt, 0pt> {41}{$\ss 4$}
\relabela <-1pt, 0pt> {42}{$\ss 4$}
\relabela <-1pt, 0pt> {43}{$\ss 4$}
\relabela <-1pt, 0pt> {51}{$\ss 5$}
\relabela <-1pt, 0pt> {52}{$\ss 5$}
\relabela <-1pt, 0pt> {53}{$\ss 5$}
\endrelabelbox}
\caption{Working through the details in the example of Figure
  \ref{figure:interesting disc}
\label{figure:working through main example}}
\end{figure}

{\bf Example} We illustrate the enumeration process for the example
which we considered earlier, in Figure \ref{figure:interesting disc}.
See Figure \ref{figure:working through main example}.

Sketches (6),(3),(1) of Figure~\ref{figure:working through main
  example} are the successive modifications of the embeddable tiled
disc which we first met in Figure \ref{figure:an example of an
  embeddable tiled surface}.  Call the three tiled discs $\cD_0,
\cD_1$ and $\cD_2$, respectively.  The initial embeddable tiled disc
$\cD_0$ has 8 vertices (at heights 1, 2, 3, 3.1, 3.2, 4, 5, and 6) and
7 singularities, labelled $a, b, c, d, e, f, g$.  Alternatively, our
enumeration process begins with the positive embeddable tiled disc
$\cD_2$ in sketch (1).  Its embedding is illustrated in sketch 2.  It
is made of discs connected by twisted bands.  The discs are
$\delta_1,\delta_2,\ldots,\delta_6$ of heights $1,2,\ldots,6$ and
radii $6,5,\ldots 1$, respectively.  Its boundary word is
$a_{6,4}a_{6,2}a_{3,2}a_{2,1}a_{5,1}^{-1}$.

In sketch (3) we have added the $ab$--tile associated to singularity
$b$ to make $\cD_1$.  We could equally well have added singularity
$f$ first, it would not matter.  The reader may find it interesting
(we did!) to try adding $f$ first, then $b$, and to see the pleasing way in
which the intermediate braids change, while both choices lead to the
same final braid (sketch (8)).

The component of $EB(\cD_1)$ that is not in $\partial D_1$ is lightly dashed
in sketch (4).  The new boundary of $\partial D_1$ (ie, after deleting the
dashed circle) is shown in sketch (5).

The reader may have wondered about the negative vertex at level 3.2.
As an aside, we now suggest that the interested reader
start with the embedded disc of sketch (2) and add yet one more disc
$\delta_{3.2}$ of $\rm{radius}=6-3.2$, at a level between $\delta_3$ and
$\delta_4$, only now with its negative side facing `up'.  The boundary
of $\delta_{3.2}$ can be seen to be homotopic through an embedded
annulus to the dashed unknot in sketch (4) of Figure
\ref{figure:working through main example}.  Let us refer to the disc
with the annulus attached as the new disc.

The new disc meets the old surface along part of its boundary---all of
the dashed unknot of sketch (4) except the dashed arc at polar angle
($b$).   Further deform the new disc slightly to take the dashed arc to
the solid arc at angle ($b$) in sketch (4).   If we now form the union of
the old surface and new disc, behold, we have a new embedded surface
whose boundary is exactly that of sketch (5)!
It's wonderful to see how the pasting can be done without introducing
any intersections with the previous surface.

In this way we obtain an embedded disc $D_1$ which is bounded by the
braid in sketch (5).  The piece of the old boundary which was on the
dashed circle and away from the wedge of 3--space near polar angle $b$
is in the interior of $D_1$, and the half-band at polar angle $b$ is
a subarc of the new boundary.

Sketches (7) and (8) show the next stage of adding another $ab$--tile
at singularity $f$ to get $\cD_0$.  Again the new unknot component is
lightly dashed.  Again there is a missing disc $\delta_{3.1}$, and it
fits right between $\delta_{3}$ and $\delta_{3.2}$ to give the
embedded surface $\cD_0$, which is exactly the surface in Figure
\ref{figure:interesting disc}, constructed in a systematic way.

The braid in sketch (8) is $\partial D_0$.   It is clearly a
4--braid, though there is some question as to how one can read off the
braid word efficiently and algorithmically because its strands skip
between 6 different discs.   There are many methods for this which
essentially involve keeping a table of which strands are in the braid
component at each singular angle.   We skip the details here.   More
information will be provided in a future paper on implementation details
of the algorithm.

The list of embeddable tiled discs which we have just constructed has
duplications.  Some redundancy is caused by non-uniqueness of
foliations; we will obtain duplicate embeddable tiled discs which
differ only by changes in foliation, which at most change the boundary
by an isotopy in the complement of the axis.  Additional redundancies
occur because non-isotopic embeddable tiled discs can have boundaries
which represent the same conjugacy class in $B_n$.  In another paper
we will consider the problem of implementing the algorithm, and at
that time we will address these issues.

\begin{remark} \rm
\label{remark:frustration}
The reader who has followed the details of the calculation which we
illustrated in Figure \ref{figure:working through main example} will
have learned that the boundary of the embeddable tiled disc which is
illustrated in sketch (6) of that Figure, and also in the bottom
sketch in Figure \ref{figure:interesting disc}, is represented by
the 4--braid:
$$a_{4,2}a_{3,2}a_{2,1}a_{4,3}a_{3,2}^2a_{2,1}
a_{3,2}^{-1}a_{4,3}^{-1}a_{2,1}^{-1}a_{3,2}^{-1}.$$
On the other hand, the closed braid which
is illustrated in the top left sketch in Figure
\ref{figure:interesting disc} is represented by the 4--braid:
$$a_{4,3}a_{3,2}^{-1}a_{2,1}^{-2}a_{3,2}^{-1}a_{3,2}a_{3,2}
a_{2,1}^2a_{3,2}a_{4,3}^{-1}a_{3,2}a_{2,1}.$$
A direct attempt to show
that these two 4--braids are conjugate in $B_4$ will very likely
lead to frustration.  Fortunately, an algorithmic solution to this
difficult problem is available.  Indeed, we have interfaced it with
our program for the algorithm.  We refer the reader to Appendix B,
where it is described, briefly.
\end{remark} 

\section{The halting theorem}
\label{section:halting} 

In the previous sections we have shown that there is an algorithm
which enumerates all foliated discs.  By reading the braid word of the
boundary of each disc we then get a list of all possible (conjugacy
classes of) unknots.  The problem which remains is to learn
when to stop testing and conclude that the knot in question is not going to
appear on the list and, hence, is truly knotted.

The boundary of a foliated disc is a closed braid, but our given knot $K$ will
in general not be a closed braid. The first step is to change $K$ to a closed
braid. There are many ways to do this. If the knot has $n$ Seifert circles and
$k_0$ crossings, the method given in Appendix~\ref{appendix:changing knots to
closed braids} is simple and it converts $K$ to a closed braid of $n$ strands and
word length
$k\leq k_0 + (n)(n-1)$. The pair $(n,k)$ is then a measure of complexity of
the given example. But $(n,v)$ is our measure of complexity of the foliated
disc. Our main task in this
section is to find a relationship between
$(n,v)$  and $(n,k)$. Clearly, $n=n$,
and we will show that there is a relationship of the form $v<f(n,k)$ where $f$
is an appropriate function.  This will prove that we need not look for
arbitrarily complicated discs.

Establishing such an upper bound on $v$ using foliated surface
techniques is an interesting and significant open problem---one we have
not solved.
We believe that a solution which is in the spirit of our algorithm exists, but
that finding it  will depend upon obtaining a better understanding than we
have at this time of the boundary words. We have been able to
prove that an upper bound exists by using the machinery in \cite{BMV}, but we
were unable to find an explicit formula.  We have not included that proof
because it requires the development of a great deal of background material.   

We present a different approach. We establish
our upper bound by first constructing a triangulation of the complement of $K$
which is adapted to closed braids. In particular, we do it so that the braid
axis meets
$\alpha$ tetrahedra, and meets them in a controlled way. We will need to count
the total number
$t$ of tetrahedra in the triangulation. After that we will use the Kneser--Haken
theory of normal surfaces to obtain an upper bound, depending on t,  on the
number of intersections of the axis with a single tetrahedron. Multiplying by
$\alpha$ we will obtain the bound that we need.

\subsection{Constructing the triangulation}
\label{subsection:constructing the triangulation}
Our goal in this section is to construct the triangulation and to prove Lemma
\ref{lemma:triangulation} below. A suitable triangulation of $S^3$ is one for
which:
\begin{enumerate}
\item The knot $K$ is in the 1--skeleton,
\item A regular neighborhood of $K$ is triangulated
\item The knot $K$ is a closed braid, and the braid axis meets a fixed
number $\alpha$ of tetrahedra, in a fixed way, ie as a straight line from
the center of one of the faces to the opposite vertex.  
\end{enumerate}
The construction of a suitable triangulation of $S^3$ involves many technical 
details. The reader
who is willing to accept the fact that we can construct one with the stated
properties will not lose the thread of the argument by proceeding at this time
directly to Lemma
\ref{lemma:triangulation}, and accepting its truth. 

In the interests of
efficiency of the algorithm,  we will try to minimize
$t$ and
$\alpha$ in our construction of a suitable triangulation. We start by defining a
few pieces of the construction.  At several points in the construction we must
subdivide.  We need several specific subdivisions to do this well.

\begin{figure}[htb!]
\let\ss\scriptstyle
\centerline{\relabelbox\small
\epsfxsize 3.5in\epsfbox{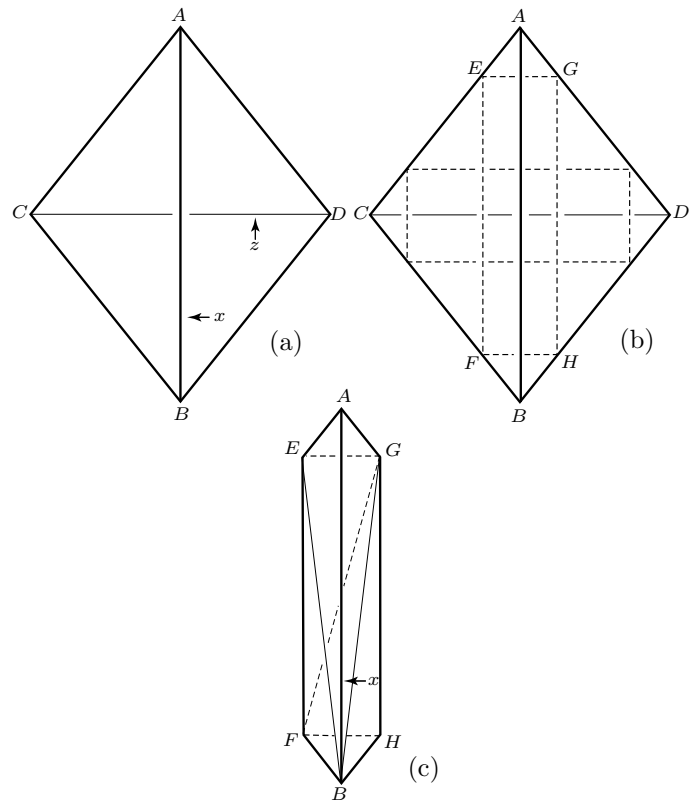}
\relabel {a}{(a)}
\relabela <5pt, -5pt> {b}{(b)}
\relabela <10pt, 0pt> {c}{(c)}
\relabela <-2pt, 0pt> {z}{$\ss z$}
\relabela <-.5pt, 0pt> {x1}{$\ss x$}
\relabela <-1.5pt, 0pt> {x2}{$\ss x$}
\relabela <-3pt, 0pt> {A1}{$\ss A$}
\relabela <-3pt, 0pt> {A2}{$\ss A$}
\relabela <-1pt, -1pt> {A3}{$\ss A$}
\relabela <-2pt, -1pt> {B1}{$\ss B$}
\relabela <-3pt, 0pt> {B2}{$\ss B$}
\relabela <-3pt, 0pt> {B3}{$\ss B$}
\relabela <-2pt, 0pt> {C1}{$\ss C$}
\relabela <-1.5pt, 0pt> {C2}{$\ss C$}
\relabela <-2.2pt, -1pt> {D1}{$\ss D$}
\relabela <-1pt, 0pt> {D2}{$\ss D$}
\relabela <-2pt, 0pt>  {E1}{$\ss E$}
\relabela <-4pt, 0pt> {E2}{$\ss E$}
\relabela <-4pt, 0pt> {F1}{$\ss F$}
\relabela <-2pt, 1pt> {F2}{$\ss F$}
\relabel {G1}{$\ss G$}
\relabel {G2}{$\ss G$}
\relabela <-1pt, 0pt> {H1}{$\ss H$}
\relabela <-1pt, 0pt>  {H2}{$\ss H$}
\endrelabelbox}
\caption{The crossing triangulation
\label{figure:crossing triangulation}}
\end{figure}

\med{\bf Definition}\qua Each crossing of the braid strands will be carried
by a tetrahedron which we subdivide with the {\em crossing
  triangulation} which is defined as follows.  Given a tetrahedron $T$
with vertices $A, B, C, D$ as in Figure~\ref{figure:crossing
  triangulation}(a), start by slicing it by two planes $P_1, P_2$
which separate the tetrahedron into 3 pieces as in
Figure~\ref{figure:crossing triangulation}(b).  The top and bottom
pieces are affine triangular prisms and contain $\overline{AB}$ and
$\overline{CD}$ respectively.  The middle piece is an affine cube
which can be further divided into two affine prisms.

An affine prism can be triangulated with 3 tetrahedra.  Consider the
prism $ABEFGH$ in Figure~\ref{figure:crossing triangulation}(c).  The
simplices $ABEG, BEFG$ and $BFGH$ are a triangulation of the prism.
After cutting the middle piece into two prisms, $T$ has been divided
into 4 prisms, each of which we can subdivide into 3 tetrahedra in a
compatible way to make a triangulation of $ABCD$ with 12 tetrahedra.
We call this subdivision the {\em crossing triangulation}.  Note that
the crossing triangulation of $T$ has 5 triangles in each original
face of $T$.

\begin{figure}[htb!]
\cl{\epsfxsize 2in\epsfbox{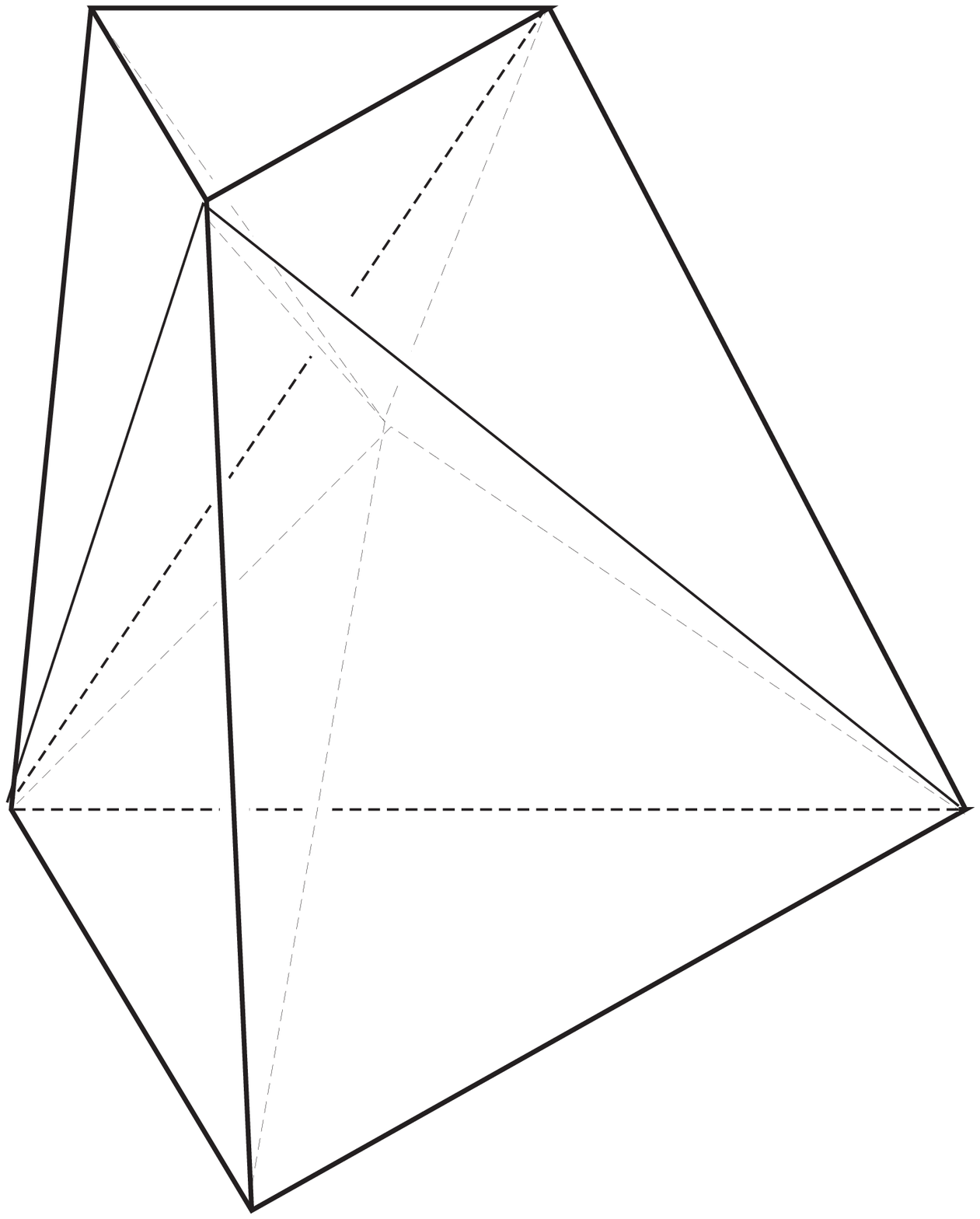}}
\caption{The frustum triangulation
\label{figure:frustum triangulation}}
\end{figure}

\med{\bf Definition}\qua Given a tetrahedron, consider the \emph{frustum}
obtained by slicing off one corner (say, the top corner) with a plane
(see Figure~\ref{figure:frustum triangulation}).  This frustum is a
convex polyhedron with 3 quadrilateral faces and two triangular faces.
The {\em frustum triangulation} is the triangulation achieved by
drawing in either of the two diagonals in each quadrilateral face,
thus triangulating the boundary, then taking the cone of the
triangulation on the boundary into any interior point.
The boundary of the triangulated frustum has 8 faces, so there are 8
tetrahedra after subdividing.  The original tetrahedron has been
subdivided into 9 tetrahedra.

\begin{figure}[htb!]
\let\ss\scriptstyle
\centerline{\relabelbox\small
\epsfxsize 4in\epsfbox{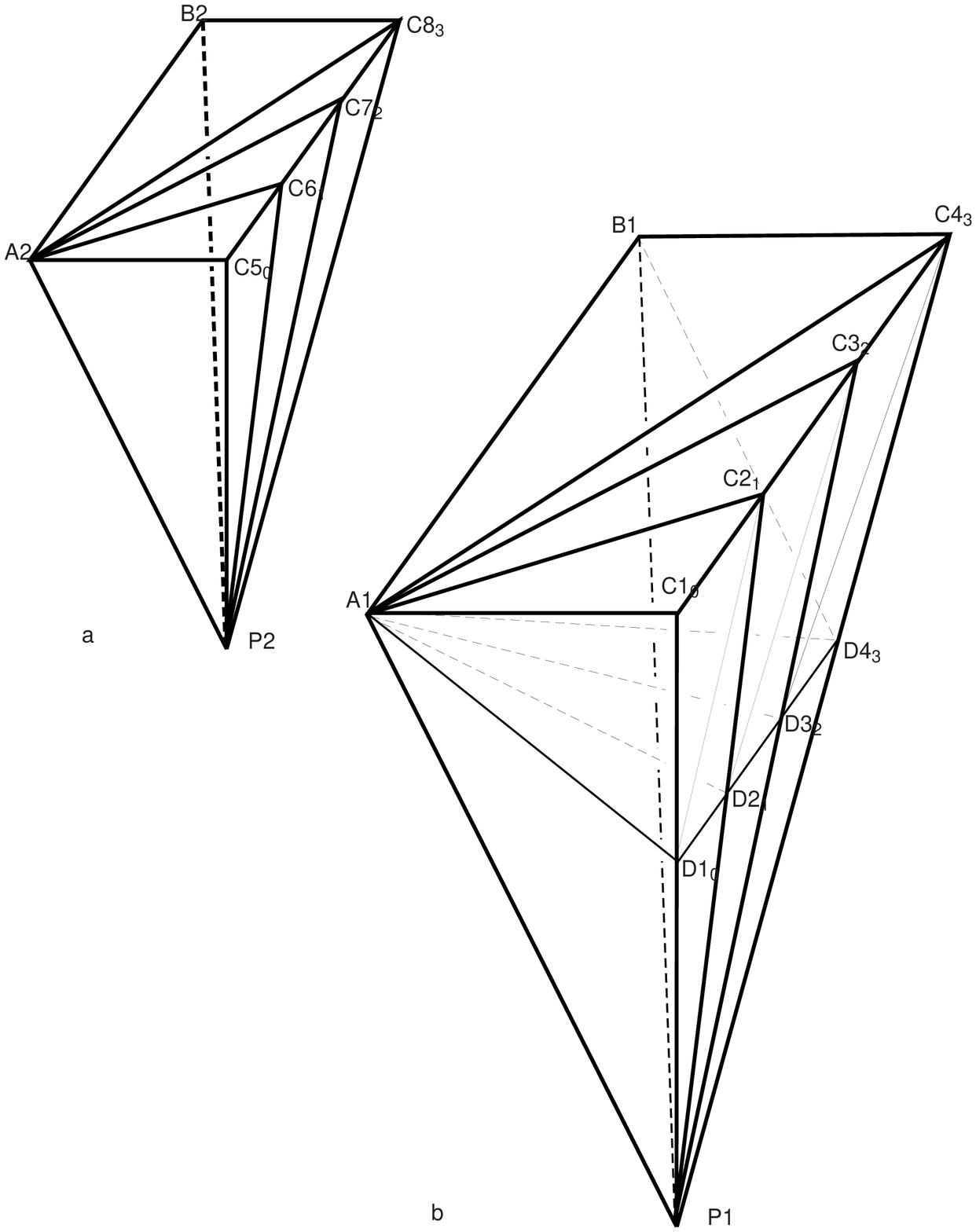}
\relabel {a}{(a)}
\relabela <5pt, 0pt> {b}{(b)}
\relabela <-3pt, 0pt> {A1}{$\ss A$}
\relabela <-3pt, 0pt> {A2}{$\ss A$}
\relabela <-2pt, -1pt> {B1}{$\ss B$}
\relabela <-3pt, 0pt> {B2}{$\ss B$}
\relabela <-4pt, -1pt> {P1}{$\ss P$}
\relabela <-3pt, 0pt> {P2}{$\ss P$}
\relabela <-2pt, 0pt> {C1}{$\ss C_0$}
\relabela <-1.5pt, 0pt> {C2}{$\ss C_1$}
\relabela <-1.5pt, 0pt> {C3}{$\ss C_2$}
\relabela <-1.5pt, 0pt> {C4}{$\ss C_3$}
\relabela <-1.5pt, 0pt> {C5}{$\ss C_0$}
\relabela <-1.5pt, 0pt> {C6}{$\ss C_1$}
\relabela <-1.5pt, 0pt> {C7}{$\ss C_2$}
\relabela <-1.5pt, 0pt> {C8}{$\ss C_3$}
\relabela <-2.2pt, -1pt> {D1}{$\ss D_0$}
\relabela <-1.5pt, 0pt> {D2}{$\ss D_1$}
\relabela <-1.5pt, 0pt> {D3}{$\ss D_2$}
\relabela <-1pt, 0pt> {D4}{$\ss D_3$}
\endrelabelbox}
\caption{The edge triangulation
\label{figure:edge triangulation}}
\end{figure}

\med{\bf Definition}\qua The final subdivision of a tetrahedron we will
define we call the {\em edge triangulation} because we will use it for
tetrahedra at the edges of our picture.  Suppose we have a rectangular
pyramid triangulated as in Figure~\ref{figure:edge triangulation}(a).
Figure~\ref{figure:edge triangulation} shows a pyramid divided into 4
tetrahedra, but in general there might be any finite number.  Pick a
point $D$ in the edge $\overline{CP}$ and consider the plane
containing points $A,B$ and $D$.  Let $D_1,\ldots, D_n$ be the points
of intersection between this plane and the edges
$\overline{C_1,P},\ldots,\overline{C_n,P}$.  This plane cuts each of
the original tetrahedra into 2 pieces, one of which (the one
containing $P$) is a tetrahedron and one which is not a tetrahedron.

The final step in the construction is to divide the objects which are
not tetrahedra by another plane.  Consider the plane containing the
triangle $AC_kD_{k-1} (0<k\leq n)$.  This cuts the object with
vertices $A,C_k, D_k, D_{k-1}, C_{k-1}$ into two tetrahedra.  The end
result is that each of the original simplices in
Figure~\ref{figure:edge triangulation}(a) is now 3 simplices in (b),
except for the tetrahedron $ABC_nP$ has been divided into 2.  Note
that
the face $ABP$ is not subdivided by this triangulation.

We are now prepared to triangulate $S^3$ so as to satisfy  conditions
1 and 2 from above.  Consider Figure~\ref{figure:triangulated
  cylinder}.  It is a triangulated rectangle, but after identification
of the top and bottom edges it becomes a triangulated cylinder or
annulus.  The bold lines represent the original knot (as a braid) and
the thinner lines are other edges in the triangulation.  The cylinder
has been draw with the middle section shown larger than the other two
only because that is where the interesting detail lies.  It should be
understood that the top and bottom portions will be larger when the
cylinder is embedded.

\begin{figure}[htb!]
\let\ss\scriptstyle
\centerline{\relabelbox\small
\epsfxsize 3.5in\epsfbox{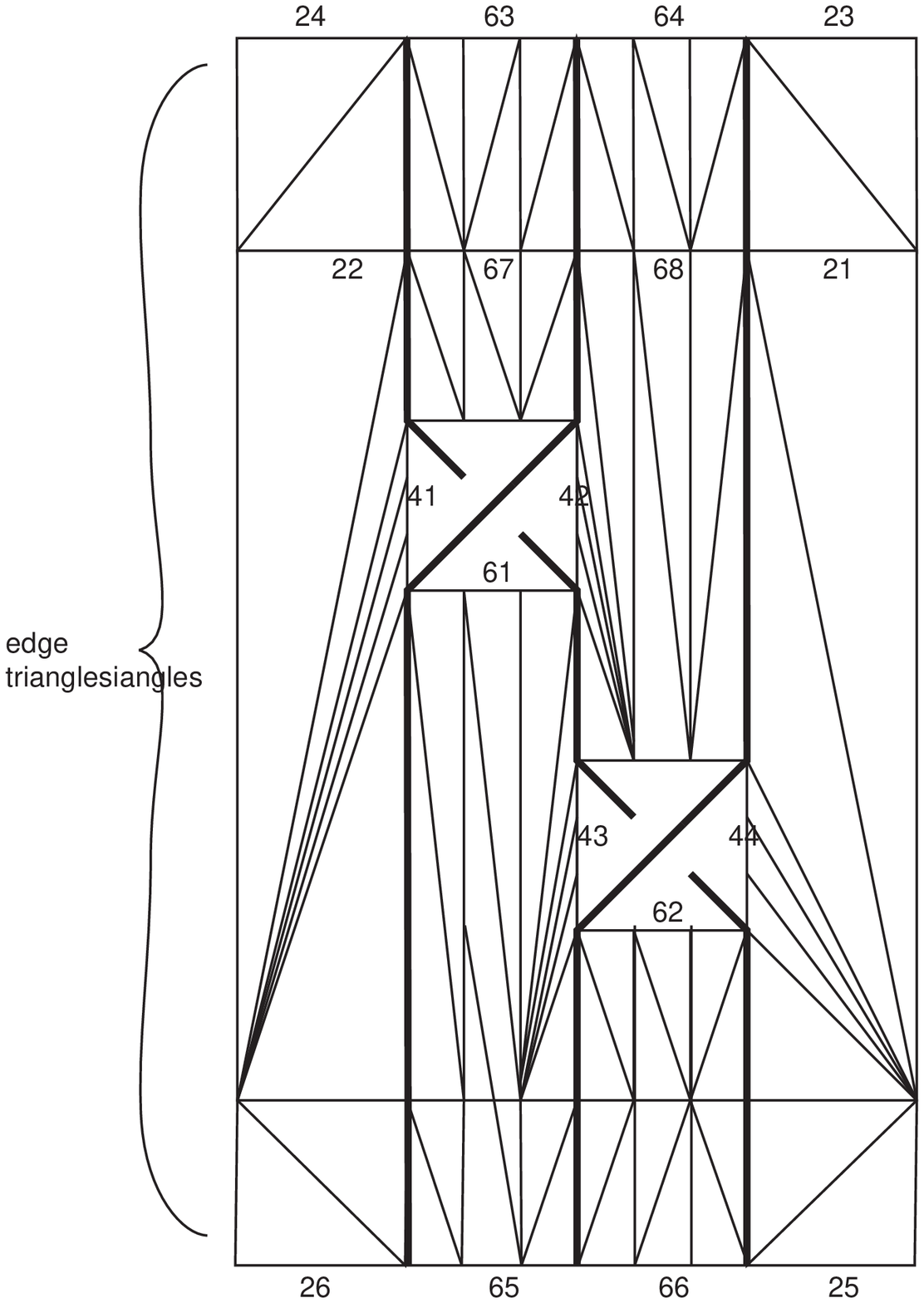}
\relabel {edge}{edge}
\relabel {triangles}{triangles}
\relabel {21}{$\ss2$}
\relabel {22}{$\ss2$}
\relabel {23}{$\ss2$}
\relabel {24}{$\ss2$}
\relabel {25}{$\ss2$}
\relabel {26}{$\ss2$}
\relabel {41}{$\ss4$}
\relabel {42}{$\ss4$}
\relabel {43}{$\ss4$}
\relabel {44}{$\ss4$}
\relabel {61}{$\ss6$}
\relabel {62}{$\ss6$}
\relabel {63}{$\ss6$}
\relabel {64}{$\ss6$}
\relabel {65}{$\ss6$}
\relabel {66}{$\ss6$}
\relabel {67}{$\ss6$}
\relabel {68}{$\ss6$}
\endrelabelbox}
\caption{The triangulated cylinder
\label{figure:triangulated cylinder}}
\end{figure}

The goal of our triangulation is to minimize $t$ and $\alpha$ while
satisfying conditions 1,2 and 3.  Satisfying 1 is easy---we build our
triangulation by starting with $K$ as out 1--skeleton and extending it.
To satisfy 2, we will use the following condition:
\paragraph{Neighborhood Condition}
\begin{itemize}
\item  \em{No two simplices which meet
    different strands of $K$ may intersect.} 
\end{itemize}
Clearly, any triangulation satisfying the Neighborhood Condition has a
triangulated regular neighborhood of $K$.

We assume that we have been presented with $K$ as a braid with $n$
strands and word length $k$ in the Artin generators.  (See Appendix~B
for an exposition on how to convert a knot to a braid efficiently.)
We think of $K$ as lying on a cylinder winding around and parallel to
the braid axis $A$, where $A$ is the $z$--axis in 3--space.
Furthermore,
think of this cylinder as having a triangular cross section, ie, it
is made of three flat sides.  All the crossings are assumed to be on
one side the cylinder, the other two sides carry the closing of the
braid into a knot.  Thus, in Figure~\ref{figure:triangulated cylinder}
all the crossings of the bold lines are in the middle portion of the
figure and in the top and bottom sections there are no crossings.

We start by considering a braid with $n$ strands and no crossings.  In
this case all three sections will be triangulated like the upper and
lower ones.  To guarantee the Neighborhood Condition, subdivide the
rectangle between each adjacent pair of strands in each section of the
cylinder into 3 rectangles, the subdivide each these into a pair of
triangles.  Thus, for each adjacent pair of strands there are 18
triangles, or $18(n-1)$ in all.

Strands 1 and $n$ are also adjacent to 3 rectangles which we divide
into 2 triangle each for 12 more triangles.  These are the \em{edge
  triangles}; the others are \em{interior triangles}.  Thus, if $K$
is an $n$ stranded braid with no crossings it lies on a triangulated
cylinder with $18n-6$ triangles.

Now consider adding a crossing.  The strands cannot cross and stay on
the cylinder, so we fatten the cylinder at each crossing by adding a
tetrahedron so that the two strands are on opposite edges of the
tetrahedron.  We triangulate the tetrahedron with the crossing
triangulation described earlier (Figure~\ref{figure:crossing
  triangulation}, thus guaranteeing the neighborhood condition.  We do
not draw this crossing triangulation in
Figure~\ref{figure:triangulated cylinder}, but note that it introduces
two new vertices on four sides of a square around the crossing.  Each
crossing breaks the rectangles between the strands into two rectangles
(plus the crossing).  We subdivide each rectangle as before,
converting each triangulated rectangle into 2, introducing 6 more
triangles.

The square around the crossing has 4 vertices on its sides, as well.
These vertices are in the middle of a vertical edge on a triangle
inside a rectangle triangulation.  We retriangulate that triangle by
adding straight lines from the new vertices to the opposite vertex of
the triangle, thus introducing 4 more triangles on each side.

The crossing triangulation has 5 triangles on each face of the
tetrahedron, and two faces are on each side of the cylinder.  Thus
each crossing introduces 6 + 2(4) +2(5) = 24 new triangles onto the
cylinder.  (In addition, each crossing triangulation has 12 tetrahedra
which we will count later.)  After adding all $k$ crossings we will
have a triangulated cylinder with $18n-6 + 24k$ triangles.  If $K$ is
a knot, strands 1 and $n$ must be in at least one crossing each,
introducing 4 additional edge triangles on each side, so at least 18
triangles are edge triangles and at most $18n + 24k -24$ are interior
triangles. 

The boundary of the triangular cylinder is a pair of triangles.  Cap
off the cylinder with two triangles to make it a sphere triangulated
with $18n + 24k -4$ triangles.  Pick two points in $A$, one enclosed
by the sphere and one in its exterior, called the South Pole (SP) and
North Pole (NP), respectively.

The cone of the triangulation on the sphere to SP is a triangulated
3--ball with $18n + 24k -4$ tetrahedra.  Unfortunately, the
neighborhood condition is now violated because all the tetrahedra meet
at a single point, SP.  We subdivide each tetrahedron coming from an
interior triangle by cutting the ``tip'' off the tetrahedron, yielding
a tetrahedron at SP and an untriangulated frustum.  We triangulate
each frustum into 8 tetrahedra by the frustum triangulation given
earlier, for a total of 9 new tetrahedron for each original one.  Thus
we get at most $9(18n + 24k -24)$ tetrahedra from interior triangles.

The tetrahedra coming from the edge triangles we triangulate with the
edge triangulation, replacing each edge tetrahedron with 2 or three
new tetrahedra.  Thus we have at most $9(18n + 24k -24)) + 3\cdot 18$
tetrahedra coming from the cylinder.  We do this edge triangulation
so that the sides of the edge tetrahedra containing the boundary of
the cylinder at not subdivided at all.

Finally, the tetrahedra made by coning the the boundary of the
cylinder do not need any subdividing.  These tetrahedra do not touch
the braid, thus they vacuously satisfy the neighborhood condition.
Furthermore, the only other tetrahedra touching one of these
tetrahedra are edge tetrahedra, and these are not subdivided on the
side touching the tetrahedra in question, so there is no need to
subdivide them at all.

Thus, our ball has $9(18n+24k-24) + 3\cdot 18 + 2$ tetrahedra.
Perform the same construction on the exterior of the sphere to yield
twice that many tetrahedra.  Now we add in all $12k$ tetrahedra in the
crossing triangulations at each crossing to get $324n + 432 k - 376$
tetrahedra in our triangulated triangulated $S^3$, which, by
construction, satisfies the neighborhood condition and hence condition
1 and 2 from the beginning of this section.  Notice that the braid axis $A$
intersects exactly 4 tetrahedra in the triangulation, and does so in a canonical
fashion.  The reason is that $A$ is the $z$--axis and passes through SP, NP and
the four tetrahedra made by coning the boundary of the cylinder. By
construction, it meets each of these tetrahedra as a line segment which runs
from the center of one of the faces to the opposite vertex.  Thus we have proved:

\begin{lemma}\label{lemma:triangulation}
  Given a closed braid $K$ of $n$ strands and word length $k$ in the
  Artin generators, there is a triangulation of $S^3$ satisfying the
  neighborhood condition with $324n + 432 k - 376$ tetrahedra, only 4
  of which meet the braid axis in their interiors. \endpf
\end{lemma}

\subsection{An upper bound on $v$}
\label{subsection:upper bound}
The next step in the formulation and proof of the halting theorem is to
obtain an upper bound for how many times the braid axis meets a single
tetrahedron in the triangulation. For this, we turn to the theory of normal
surfaces in 3--manifolds. A normal surface in a
triangulated 3--manifold is a PL surface in the complement of the
0--skeleton with some restrictions on how it can intersect a
tetrahedron of the triangulation. 

See \cite{Hass} for a review of the essential facts
which we need about normal surfaces.  Among them is the fact that the disc which
our unknot bounds can be isotoped to a normal surface. Also, each  component of
the intersection of that normal surface with a tetrahedron is a planar triangle
or quadrilateral with corners in the interiors of the distinct edges of the
tetrahedron.  Also, for each tetrahedron there are 7 different combinatorial
types of possible intersections.    

\begin{lemma}
\label{lemma:HLP bound} Let
$M$ be the closure of
$S^3$ minus the triangulated regular neighborhood of an unknot $K$.
Let $t$ be the number of 3--simplices in the triangulation of $M$ and
let $S$ be any simplex in the triangulation of $M$.  Then some longitude of
$K$ on $\partial M$ bounds a disc $F$ which is a normal surface, and 
the number
of components of $S\cap F$ is bounded above by $7t 2^{7t+2}$.
\end{lemma}

\pf Fortunately 
we do not need to do any work at all. Lemma 6.1 of \cite{HLP} uses a very
different triangulation of a knot complement from ours, but provides
exactly the estimate we need: that the number of components of a
given combinatorial type in $S\cap F$ is at most $t2^{7t+2}$.  
We thank Joel Hass for several useful discussions on this matter. Since there
are 7 possible combinatorial types the assertion follows.  
\endpf

Our Halting Theorem is an immediate corollary of Lemmas
\ref{lemma:triangulation} and \ref{lemma:HLP bound}. It says that we
can stop looking for tiled discs if one hasn't shown up before a given time.

\begin{hthm}\label{theorem:upper bound}
  Given an unknotted closed braid $K$ of $n$ strands and word length
  $k$ in the Artin generators, let $t=324n + 432 k - 376$.  Then there is an
embedded disc in $S^3$ whose boundary
  is a closed braid conjugate to $K$ and the induced braid foliation
  has complexity no higher than $(n, \ 28 t 2^{7t+2})$.
\end{hthm}

\pf  By the arguments which are
reviewed in Section
\ref{section:braid
 foliations of spanning surfaces for knots} of this paper, the normal
surface can be isotoped to a tiled surface $\cD=(D,G,C),$  via an
isotopy which takes place outside a neighborhood of $A$. (If $\cD$ is
inessential, there is an essential tiled disc
  with fewer intersections with the braid axis.)  
As the isotopy is away from the braid axis, the points of intersection
of the surface with $A$ are unchanged.    

By Lemma \ref{lemma:triangulation}, the triangulation has at most
$t=324n + 432 k - 376$ tetrahedra. By Lemma~\ref{lemma:HLP bound} the
normal surface may be assumed to meet each of the $\leq t$ tetrahedron
in at most $7t 2^{7t+2}$ components.  By Lemma
\ref{lemma:triangulation} the braid axis meets 4 tetrahedra. By our
construction of the triangulation in Section
\ref{subsection:constructing the triangulation} the braid axis $A$
passes through each of these as a line segment which meets each sheet
of the normal surface at most once, transversally.  The total number
of intersections of $A$ with $\cD$ is then bounded above by $28 t
2^{7t+2}$. Thus the complexity $(n,v)$ of the disc we seek is at most
$(n, \ 28 t 2^{7t+2})$, where $t\leq 324n + 432 k - 376$. \endpf

\section{The algorithm}
\label{section:the algorithm}

We summarize our test for whether $K$ {\em is} the unknot in the form
of an algorithm which is based on our work in earlier sections.  We
are given a knot $K$, which we may assume is given as a braid with $n$
strands (see Appendix~\ref{appendix:changing knots to closed braids})
and Artin word length $k$.  As in the previous section, let $t = 324n
+ 432 k - 376$, so there is a triangulation of the knot complement
with
fewer than $t$ tetrahedra. 

Let $[K]$ denote the conjugacy class of the braid word given by $K$.
The algorithm for determining whether $K$ is the unknot is outlined
below in Table~\ref{figure:the algorithm} in pseudo-code.

\begin{theorem}
The algorithm in Table~\ref{figure:the algorithm} will enumerate a
list of closed braids containing 
examples from each conjugacy class of closed braids representing the
unknot.   The given knot is unknotted if and only if it is in the 
conjugacy class of one of
the braids on our list.   
\end{theorem}

\begin{table}\small\cl{%
\fbox{\parbox{0.98\textwidth }{
\noindent{\tt
(1) for each integer $P$ from $n$ to ${14t 2^{7t+2}} + n$\\
(2){\mbox \qua} for each positive embeddable tiled disc $\cD'$ \\
\hspace*{1 in}  with exactly $P$ positive vertices\\
(3){\mbox \qua} {\mbox \qua} for each way of adding $P-n$ negative
vertices to 
$\cD'$ to get $\cD$\\
(4){\mbox \qua} {\mbox \qua} {\mbox \qua} if $\cD$ is embeddable
and
all b-arcs are essential\\
(5){\mbox \qua} {\mbox \qua} {\mbox \qua} {\mbox \qua} find the braid word $W(\partial \cD)$\\
(6){\mbox \qua} {\mbox \qua} {\mbox \qua} {\mbox \qua} compute its conjugacy class $[W(\partial\cD)]$\\
(7){\mbox \qua} {\mbox \qua} {\mbox \qua} {\mbox \qua} if  $[W(\partial\cD)] = [K]$\\
(8){\mbox \qua} {\mbox \qua} {\mbox \qua} {\mbox \qua} {\mbox \qua} $K$ is an unknot\\
(9){\mbox \qua} {\mbox \qua} {\mbox \qua} {\mbox \qua} endif\\
(10){\mbox \qua} {\mbox \qua} {\mbox \qua} endif\\
(11){\mbox \qua} {\mbox \qua} next $\cD$\\
(12){\mbox \qua} next embeddable tiled disc\\
(13) next $P$
}
}}}
\caption{The algorithm to test for unknottedness
 \label{figure:the algorithm}}
\end{table}

\pf
If $K$ is unknotted, then the algorithm is guaranteed to find a disc
that $K$ bounds.  Theorem~\ref{theorem:constructing all embeddable
  tiled discs} proves that lines (1), (2), (3) will generate all tiled
discs which are embeddable tiled discs.
Theorem~\ref{theorem:embeddability test} and
Proposition~\ref{proposition:inessential b-arcs} give the algorithmic tests for
embeddability and essential $b$--arcs used in line (4).
Theorem~\ref{theorem:boundary word} explains how to read the braid
word as needed in line (5).  The solution to the conjugacy problem
needed for lines (6) and (7) can be found in \cite{BKL}.  The Halting
Theorem (Theorem~\ref{theorem:upper bound}) tells us that if $K$ does
bound a tiled disc, then that disc has no more than $28t 2^{7t+2}$
vertices.  Since it must have $n$ more positive than negative
vertices, we only need to check up to ${14t 2^{7t+2}} + n$
positive vertices.
\endpf

\begin{remark}\rm
\label{remark:generating examples}
The algorithm in Table~\ref{figure:the algorithm} can easily be
modified to generate 
examples.  One could generate all unknotted braids with a certain size
embeddable tiled disc, for example.  In this way one could build some
{\em unknot} tables which might speed up the work of others who are
building {\em knot} tables.  One could search for unknots whose
embeddable tiled disc satisfies a particular property, such as not
having any valence one vertices in the graph.  Such examples are
interesting in that they are the unknotted closed braids which have no
trivial loop.  One could search for any property of the embeddable
tiled disc (including all properties of the braid word itself).  We
will discuss this more fully in our paper on the implementation of the
algorithm.
\end{remark}

\med{\bf Running Time and Implementation Issues}\qua  This program has
been implemented and is running, but not yet robustly tested.  There
are many aspects of an implementation that this paper has ignored,
such as how to represent embeddable tiled discs on a computer, how to
enumerate positive embeddable tiled discs, and how to decide where to
add negative vertices.  There are many solutions to these problems,
and our current implementation is crude at best.  When the entire
program is in more mature form we will write a paper discussing these
issues and giving experimental results.

Typically, one would like an estimate of running time for an
algorithm.  In this case it is best not contemplated.  There are
$2^{P-2}$ labelled trees (representing the positive embeddable tiled
discs of the first step of the algorithm with cyclically ordered
vertices) each of which has $2^{P-1}$ possible signs and $(P-2)!$
different cyclic orders on the singularities.  Thus, the first part of
the algorithm (generating all positive embeddable tiled discs of $P$
vertices)
has running time at of at least $P^{P-2}!(P-2)!2^{P-1}$, which assumes
that everything is generated with perfect efficiency and no redundancy
and doesn't take into account multiple, non-isomorphic embeddings of
the graphs in the disc.

Adding the negative vertices is no faster, just harder to analyze.
The good news is that there are many more non-embeddable or
inessential ways to add negative vertices than embeddable, essential
ways.  Thus the number of embeddable tiled discs starts going down as
negative vertices are added.  The bad news is that the algorithm is
still exponential time.

In practice, we don't bother with implementing the upper bound in line
(1).  Because of the exponential growth in the number of examples,
this algorithm implemented on a single computer won't be able to go
beyond, say, 20 vertices in any reasonable time.

The test for embeddability is quadratic in the number of vertices if
implemented as written, as is the method for extracting the braid
word.  A faster version (in progress) of the embeddability test would
be very helpful.  The algorithm also depends heavily on the conjugacy
computation of \cite{BKL}, and the running time of that computation
has not yet been fully analyzed.  The running time of the related
solution to the word problem in \cite{BKL} has been fully analysed and
is shown in \cite{BKL} to be ${\cal O}(|K|^2n)$.

\newpage

\begin{appendix}\small\parskip 4pt plus 2pt minus 2pt

\makeatletter
\let\@@itemize@\itemize
\def\itemize{\@@itemize@\parskip 0pt\relax}
\let\@@enumerate@\enumerate
\def\enumerate{\@@enumerate@\parskip 0pt\relax}
\def\@listi{\leftmargin25.5pt\parsep 0pt\topsep 0pt 
 \itemsep2pt plus2pt minus1pt}
\let\@listI\@listi
\@listi
\makeatother

\section{Appendix: \ Changing knots to closed braids}
\label{appendix:changing knots to closed braids}

 In spite of a  flood of
recent applications of braid theory to knot theory, many topologists 
regard the problem of changing knot diagrams to closed braids as a
difficult project.  We give, here, a very simple algorithm which does
the job.  It is due to Vogel \cite{Vogel}, and incorporates
earlier ideas of Yamada \cite{Yamada}.  

Let $D(K)$ be a $c$--crossing diagram on the plane $\reals^2$ which
describes a knot $K$.  Smoothing the $c$ crossings, the diagram is
replaced by a {\it Seifert diagram}, ie  by a collection of oriented
Seifert circles $s_1,\dots,s_n$ which are joined in pairs by
half-twisted bands.  We indicate the bands symbolically by signed {\it
  ties}, where the sign depends on the sign of the associated
crossing.  An example is given in Figure \ref{figure:Seifert diagram
  1}.  (In the examples we omit the signs on the ties because they are
irrelevant to the present discussion.)

\begin{figure}[htb!]
\let\ss\scriptstyle
\centerline{\relabelbox\small
\epsfxsize 4in\epsfbox{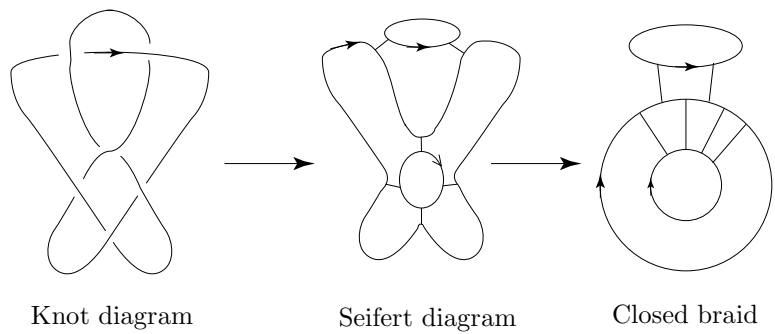}
\relabela <-10pt,0pt>  {Knot}{Knot diagram}
\relabela <-10pt,0pt>  {Seif}{Seifert diagram}
\relabela <-10pt,3pt> {Closed}{Closed braid}
\endrelabelbox}
\caption{A Seifert diagram
\label{figure:Seifert diagram 1}}
\end{figure}

Each pair of Seifert circles $s_i,s_j$ cobounds a unique annulus
$A_{i,j}$ on the 2--sphere $\reals^2\cup\{\infty\}$.  Seifert circles
$s_i$ and $s_j$ are {\em coherently oriented} if they represent the
same element in $H_1(A_{i,j};\ints)$, otherwise they are incoherently
oriented.  A key fact, first noted by Yamada \cite{Yamada} is that if
all pairs of Seifert circles are coherent, then the diagram is already
essentially a closed braid diagram.  This is the case in the example
in Figure \ref{figure:Seifert diagram 1}.  Every pair of Seifert
circles is coherently oriented, and there is no obstructions to
sliding the ties together and grouping them into a single block.
After so-doing we can number the strands and find a representing braid
word.  To be sure, it may be necessary to pull some of the Seifert
circles through the point at infinity to obtain a more conventional
closed braid diagram, as illustrated, but the braid word can be found
without that final step, so the key point is to convert the given
diagram to one in which every pair of Seifert circles is coherently
oriented.

\begin{figure}[htb!]
\let\ss\scriptstyle
\centerline{\relabelbox\small
\epsfxsize 3.5in\epsfbox{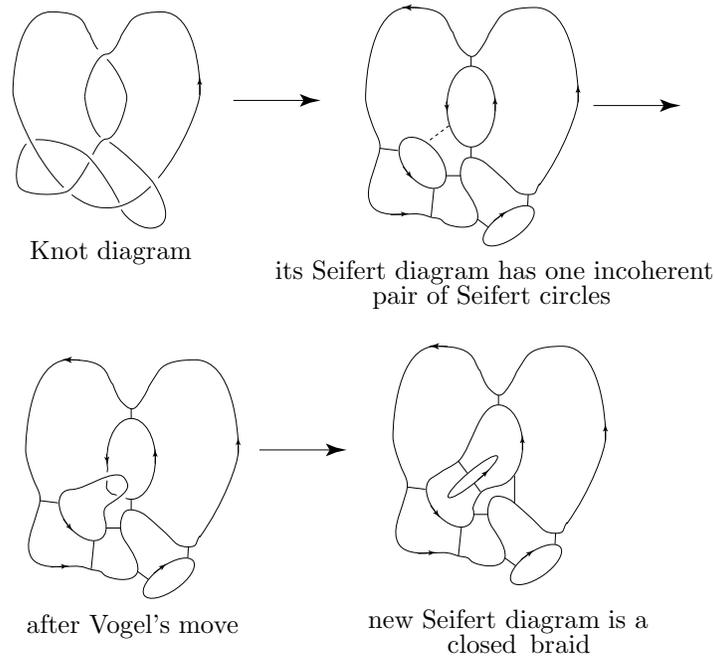}
\relabela <-5pt,0pt>  {knot}{Knot diagram}
\relabela <-10pt,0pt>  {its}{its Seifert diagram has one incoherent}
\relabela <-20pt,0pt> {of}{pair of Seifert circles}
\relabela <-5pt,0pt> {After}{after Vogel's move}
\relabela <-10pt,0pt> {New}{new Seifert diagram is a}
\relabela <20pt,0pt> {closed}{closed \,braid}
\endrelabelbox}
\caption{Another example
\label{figure:Seifert diagram 2}}
\end{figure}

In the second example, given in Figure \ref {figure:Seifert diagram 2}
  there is an incoherent pair of Seifert circles.  

 Yamada proves
that if an incoherent pair of Seifert circles exists, then
one can always find a pair, such as the pair which is joined by
the dotted arc in the second sketch for example (b).  Call the
circles $s_p,s_q$ and call the arc $\alpha_{p,q}$.  We require
that the interior of
$\alpha_{p,q}$ be disjoint from $D(K)$.    Vogel's contribution
was to notice that we can then use $\alpha_{p,q}$ to define a
Reidemeister move of type II, as in the third sketch.  This
move preserves the number of Seifert circles and reduces the
number of incoherent pairs.  Thus after a finite
number of Vogel moves we will obtain a new diagram
$D^\prime(K)$ in which every pair of Seifert circles is coherently
oriented.

Notice that the number $n$ of Seifert circles is unchanged by the
Vogel moves, and when all the Seifert circles are coherently
oriented this number is the braid index of the given example.  Hence,
we may compute $n$ from the given knot
diagram of $K$.  The writhe of the diagram is also unchanged by
Vogel moves, and we finally obtain a closed braid the writhe is
the number we have called the exponent sum.  
Thus, given a diagram of $K$ we can easily convert it to a braid for
use in our algorithm.   The crossing number, is increased by Vogel
moves, but it is proved in \cite{Vogel} that the increase is at most
$(n-1)(n-2)$.

\section{Appendix: \ Testing for conjugacy}
\label{appendix:testing for conjugacy}

In this appendix we describe, briefly, the solution to the conjugacy
problem which is given in \cite{BKL}.  We choose that solution over
the other known solutions, for example \cite{ElM}, because it uses the
band generators and so is natural for our work.  Also, there is hope
that when various technical difficulties are overcome it can be proved
to be a polynomial-time algorithm.

Let $a_{t,s}$ be one of the band generators of the braid group which
were introduced in Section 3 of this paper and illustrated in Figure
\ref{figure:band generators}.  Notice that these generators include
the more familiar generators $\{ \sigma_1,
\sigma_2,\dots,\sigma_{n-1}\}$ as a proper subset, because $a_{i+1,i}
= \sigma_i$.  It is proved in \cite{BKL} that:

\begin{prop}
\label{prop:presentation}
$B_n$ has a presentation with generators $\{a_{ts} ; \ n\geq t>s\geq
1 \}$ and with defining relations:
\begin{enumerate}

\item \label{eqn:bandrelation1}
$a_{ts}a_{rq}=a_{rq}a_{ts}\quad\hbox{\sl if}\quad (t-r)(t-q)(s-r)(s-q)>0$

\item \label{eqn:bandrelation2}
$a_{ts}a_{sr}=a_{tr}a_{ts}=a_{sr}a_{tr}\quad\hbox{\sl for\ all}\quad t,s,r
\quad\hbox{\sl with}\quad n\geq t>s>r \geq 1$.    

\end{enumerate}
\end{prop}

The defining relations in the presentation of
Theorem~\ref{prop:presentation} for $B_n$, like the more standard
presentation of Artin, only involves positive powers of the
generators.   This allows us to introduce a semigroup $B_n^+$ which
has the same presentation as $B_n$ in terms of the band generators.
There is a natural map $B_n^+\to B_n$ which takes each generator
$a_{t,s}$ of $B_n^+$ to the corresponding generator of $B_n$.
Following methods pioneered by Garside, the following embedding
theorem is proved in
\cite{BKL}.   We remark that the theorem does not appear to follow
directly from the earlier work of Garside.

\begin{prop} 
\label{proposition:embedding}
If two positive words in the band generators represent the same 
element of $B_n$,
then they also represent the same element of $B_n^+$.  
\end{prop}

In view of Proposition \ref{proposition:embedding}, we may regard $B_n^+$ as a
submonoid of $B_n$.   Generalizing the ideas of Garside (but many of the
details are different because of the new generating set) the
{\it fundamental word}
$\delta$ is introduced in
\cite{BKL}:    
\begin{equation}
\label{eqn:delta}
 \delta  = a_{n(n-1)}a_{(n-1)(n-2)}\cdots a_{21}.  
\end{equation}  
The reader who is familiar with the mathematics of braids will recognize 
that $\delta^n$ generates the center of $B_n$.    
The following two properties of $\delta$ are established in
\cite{BKL}.

\begin{lemma}
\label{lemm:delta}  
Let $\delta$ be the fundamental braid.   Then:
\begin{enumerate}
\item [{\rm (1)}] For each generator $a_{t,s}$ of $B_n$ there exists a positive
word $P_{t,s}$ in the band generators such that $a_{t,s}^{-1} =
\delta^{-1}P_{t,s}$.   
 \item [{\rm (2)}] The braid $\delta$ has the following (weak) commutativity
properties:
$a_{t,s}\delta =
\delta a_{t+1,s+1}$ and
$a_{t,s}\delta^{-1} = \delta^{-1}a_{t-1,s-1}$ 
\end{enumerate} 
\end{lemma} 

 Using (1) of Lemma \ref{lemm:delta}, we may replace all negative letters
in any word $W$ in the band generators by positive words and powers of
$\delta^{-1}$.   Using (2) we may then move all powers of $\delta^{-1}$
to the left.   In this way, one may choose an arbitrary word $W$ in the
band generators and change it to an equivalent word of the form
$\delta^pP$, where $P\in B_n^+$.

An element of $B_n^+$ is a {\it descending cycle} if it can be
represented by a word of the form
$$a_{(t_k,t_{k-1})} a_{(t_{k-1},t_{k-2})},\dots, a_{(t_2,t_1}), \
\ n\geq t_k>t_{k-1}>\dots >t_1\geq 1.$$ 
We describe a descending cycle by the array of subscripts
$(t_k,t_{k-1},\dots,t_1)$.   Two descending cycles
$(t_k,t_{k-1},\dots,t_1)$ and $(s_m,s_{m-1},\dots,s_1)$ are {\it
parallel} if no pair $(t_j,t_i)$ separates any pair $(s_q,s_p)$.
These concepts are investigated in detail in \cite{BKL}.   A {\it
canonical factor} is an element of $B_n^+$ which can be represented by
a product of parallel descending cycles.   It is proved in \cite{BKL}
that the canonical factors are precisely the `initial segments' of
$\delta$, a concept which will be familiar to those readers who have
worked through the details of the papers on Garside's algorithm.   The main
result in \cite{BKL} is:

\begin{theorem}{\rm\cite{BKL}} 
\label{theorem:conjugacy problem}
Let $W$ be an arbitrary word in the band generators which represents
$\cW\in B_n$.   Then there is a algorithmic procedure which allows one
to find, starting with $W$,  all possible representatives of the conjugacy class $[{\cal
W}]$ of the form $\delta^p A_1A_2\dots A_k,$ such that:
\begin{enumerate}
\item $p$ is maximal for all such representations,
\end{enumerate} 
and simultaneously
\begin{enumerate}
\item[\rm(2)]$k$ is minimal for all such representations. 
\end{enumerate}
Also:
\begin{enumerate} 
\item [\rm(3)] each $A_i\in B_n^+$ is a canonical factor, 
\item [\rm(4)] each product $A_iA_{i+1}, \ i=1,\dots,k-1$ is left-weighted, 
as defined in \cite{BKL}.  
\end{enumerate}
The set of finitely many words in the form $\delta^p A_1A_2\dots A_k$ 
which represent $\cW$ is known as the super summit set of
$[\cW]$.   Two elements ${\cal W}, {\cal W}^\prime
\in B_n$ are conjugate if and only if they have the same values of $p$ and $k$ and the
same super summit set.  
\end{theorem}
\end{appendix}  

\end{document}